\documentclass[12pt, a4paper]{amsart}
\allowdisplaybreaks[1]
\usepackage[all]{xy}
\usepackage{bbm}
\usepackage{mathrsfs}
\usepackage{amsmath}
\usepackage{amsfonts}  
\usepackage{amscd}
\usepackage{amssymb}
\usepackage{verbatim}%
\usepackage{hyperref}
\hypersetup{colorlinks=true,
linkcolor=blue, citecolor=red} 
\usepackage{mathtools}
\usepackage{color}
\usepackage{tensor}
\usepackage{cancel}
\usepackage{bm}

\newtheorem{thm}{Theorem}[section]
\newtheorem{lem}[thm]{Lemma}
\newtheorem{coro}[thm]{Corollary}
\newtheorem{prop}[thm]{Proposition}

\newtheorem{rem}[thm]{Remark}

\numberwithin{equation}{thm}%

\def\ggg{\mathfrak{g}}

\def\ppp{\mathfrak{p}}
\def\nnn{\mathfrak{n}}
\def\hhh{\mathfrak{h}}
\def\bbb{\mathfrak{b}}
\def\uuu{\mathfrak{u}}
\def\mmm{\mathfrak{m}}
\def\ker{{\rm Ker\,}}
\def\im{{\rm Im\,}}

\def\Hom{{\rm Hom}}
\def\End{{\rm End}}
\def\Ext{{\rm Ext}}
\def\id{{\textrm{id}}}
\def\Soc{{\rm Soc}}
\def\Rad{{\rm Rad}}

\def\bbf{\mathbb F}
\def\bbz{\mathbb Z}

\def\bbd{\mathbb D}
\def\bbv{\mathbb V}

\def\mod{\mathrm{mod\;}}
\def\Ca{{\mathcal{C}_A}}
\def\rank{\mathrm{rank}}

	\title[]{Modular representation of Reductive Lie algebras and related combinatorial category}
	\author{An Zhang}
	\thanks{This work is supported by the National Natural Science Foundation of China (NSFC) Gr:  12071136.}
	\address{School of Mathematical Sciences, East China Normal University, Shanghai 200241, China}
	\email{52215500003@stu.ecnu.edu.cn}
\begin{document}	
\maketitle
	\begin{abstract}
	We introduce and study a ``combinatorial" category related to the representations of reduced enveloping algebras of reductive Lie algebras in ``standard Levi form". It is compatible with the
	so-called AJS category in \cite{AJS94}, where AJS category is an important role in
	 studying the Lusztig's conjecture on characters of irreducible modules of algebraic group.
	\end{abstract}
\tableofcontents

\section{Introduction}
Let $G$ be a connected reductive algebraic group over an algebraically closed field $k$ of characteristic $p$. Let $\ggg=\mathrm{Lie\;}G$ be the corresponding Lie algebra, which is naturally a restricted Lie algebra. The $p$-character of simple $\ggg$-module is a linear functin $\chi$ on $\ggg$ such that $x^p-x^{[p]}-\chi(x)^p$ annihilate the simple module for all $x\in\ggg$. Then $p$-characters help us reduce
the representations of $\ggg$ to the reduced enveloping algebra $U_\chi(\ggg)=U(\ggg)/(x^p-x^{[p]}-\chi(x)^p\mid x\in\ggg)$. To figure out the representation of $\ggg$, it is enough to study the representation of each single reduced enveloping algebra. The general results of Kac and Weisfeiler (\cite{KW71}) continue to reduce the representation of reduced enveloping algebra into the case $p$-character is nilpotent, where nilpotent $p$-character means it vanishes on some Borel subalgebra of $\ggg$. With our ``combinatorial" category, we obtain a description of the character of simple modules in ``standard Levi form"

The restricted enveloping algebra is a special case of reduced enveloping algebra such that $p$-character is zero. The representations of restricted enveloping algebras are closely related to the rational representation of $G$. In particular, Lusztig's conjecture on the algebraic group with expression of the characters of the simple modules (\cite{Lu80}) has an equivalent description in restricted enveloping algebra. 
For the general nonzero $p$-character, there is a certain case that has a precise classification of simple modules by the highest weights. Following Friedlander and Parshall \cite{FP90}, this case is called the ``standard Levi form". The restricted enveloping algebra is also a special case of the ``standard Levi form". Lusztig \cite{Lu97} also proposed a conjecture of the character of simple modules in ``standard Levi form". 

To attach Lusztig'conjecture on the irreducible character over the reductive algebraic groups as well as restricted enveloping algebras, Anderson, Jantzen, and Soergel \cite{AJS94} developed some results in a graded category of restricted enveloping algebra over some commutative rings. They relate the multiplicities of simple modules with the rank of some object in a combinatorial category, which has been called the Anderson-Jantzen-Soergel category. Recently, Abe \cite{Abe} constructed the Hecke actions on the Anderson-Jantzen-Soergel category and related the multiplicities of simple modules with coefficients in the Soergel bimodule.

The graded category of reduced enveloping algebra in ``standard Levi form" is investigated by Jantzen \cite{Ja98, Ja00, Ja04} and Westaway \cite{We22, We23}. 
They give some restrictions on the structural map in their results. We follow their constructions and develop a combinatorial category compatible with the Anderson-Jantzen-Soergel category, which means it is equivalent to the Anderson-Jantzen-Soergel category when $p$-character is zero.
 We hope there is also a Hecke action on it and solve Lusztig's conjecture on ``standard Levi form" in terms of coefficients in Soergel bimodules.


\section{The graded categories}

\subsection{Preliminaries} 
  
Throughout this paper, we always assume $k$ is an algebraically closed field of prime characteristic $p$. Let $G$ be a connected  algebraic group over $k$ and set $\ggg=\mathrm{Lie\;}G$ satisfying the hypotheses$\colon$
\begin{enumerate}
\item The derived group of $G$ is simply connected.
\item The prime $p$ is good for $G$.
\item There exists a non-degenerate $G$-invariant bilinear form of $\ggg$.
\end{enumerate}
 Then $\ggg$ has a structure as a restricted Lie algebra over $k$.
  Fix a maximal toral $T$ in $G$ and set $\hhh=\mathrm{Lie\;}T$. Let $X=X(T)$ be the additive group of characters and $R\subset X$ be the root system of $G$ with respect to $T$. 
  Denote $\ggg_\alpha$ the corresponding root subspace of $\ggg$ for $\alpha\in R$. 
  Choose a positive system $R^+$, we have a triangular decomposition $\ggg=\nnn\oplus\hhh\oplus\nnn^-$, where $\nnn=\oplus_{\alpha\in R^+}\ggg_\alpha$ and $\nnn^-=\oplus_{\alpha\in R^+}\ggg_{-\alpha}$. 
  Set $\bbb=\nnn\oplus\hhh$ be the corresponding Borel subalgebra.
  
  For each $\alpha\in R$, let $\alpha^\vee$ denote the corresponding coroot. The affine reflection on $X$ is given by $s_{\alpha,mp}(\lambda)=\lambda-(\langle\lambda,\alpha^\vee\rangle-mp)\alpha$. Write $W_p$ the affine Weyl group of $R$ generated by all $s_{\alpha,mp}$ with $\alpha\in R$ and $m\in\bbz$, and write $W$ the finte Weyl group generated by $s_\alpha=s_{\alpha,0}$ with $\alpha\in R$. Denote the longest element in $W$ by $w_0$. The dot action of $W_p$ on $X$ is given by $w\cdot\lambda=w(\lambda+\rho)-\rho$ 
  , where $\rho$ is half the sum of all positive roots.
    
 Each $\lambda\in X$ defines a linear form $d\lambda$ on $\hhh$ by taking the derivative, which induces an embedding from $X/pX$ into $\hhh^*$. For each $\alpha\in R$, we set $h_\alpha$ in $\hhh$ satisfy $\langle\lambda,\alpha^\vee\rangle\equiv d\lambda(h_\alpha) (\mod p)$. Choose $x_\alpha\in\ggg_\alpha$ for each $\alpha\in R$ such that $[x_\alpha,x_{-\alpha}]=h_\alpha$.
    
 We denote by $U(\ggg)$ be the universal enveloping algebra of $\ggg$. For each $\chi\in\ggg^*$, the reduced enveloping algebra $U_\chi(\ggg)$ is the  quotient of $U(\ggg)$ by the ideal $I_\chi$ generated by all $x^p-x^{[p]}-\chi(x)^p$ with $x\in\ggg$. Each simple $\ggg$-module is in a natural way a $U_\chi(\ggg)$-module for some $\chi\in\ggg^*$. Let $S(\hhh)=S$ be the symmetric algebra over $\hhh$, which is isomorphic to $U(\hhh)$ in a natural way.
 
 We call a linear function $\chi\in\ggg^*$ is of standard Levi form if $\chi(\bbb)=0$ and there exists a subset $I$ of all simple roots such that 
 $$\chi(x_{-\alpha})=\begin{cases}
 1, & \text{ if } \alpha\in I,\\
 0, & \text{ if } \alpha\notin I.
 \end{cases}
 $$
 We fix  $\chi\in\ggg^*$ of standard Levi form and $I$ as in the definition. 
 Note that $R_I=R\cap\bbz I$ is also a root system and we write $R^+_I=R^+\cap\bbz I$. We have a Levi subalgebra $\ggg_I=\nnn_I\oplus\hhh\oplus\nnn^-_I$, where $\nnn_I=\oplus_{\alpha\in R^+_I}\ggg_\alpha$ and $\nnn^-_I=\oplus_{\alpha\in R^+_I}\ggg_{-\alpha}$. The affine Weyl group $W_{I,p}$ of $R_I$ is generated by $s_{\alpha,mp}$ with $\alpha\in I$ and $m\in\bbz$. The finite Weyl group $W_I$ is generated by $s_{\alpha,0}$ with $\alpha\in R_I$ and we denote by $w_I$ the longest element in $W_I$.
 
 Set $\ppp=\nnn_I^-\oplus\hhh\oplus\nnn$ be the parabolic subalgebra in $\ggg$, we have $\bbb\subset\ppp$.

\subsection{Definition of graded categories}
Let $J$ be the ideal of $U(\ggg)$ generated by $X^p-X^{[p]}-\chi(X)^p$ with $X\in\nnn\oplus\nnn^-$. Set $U_\chi=U_\chi(\ggg)/J$. Let $U^0$ be the image of $U_\chi(\hhh)$ in $U_\chi$, which is isomorphic to $U(\hhh)$. The PBW theorem shows that 
$$U_\chi\simeq U_\chi(\nnn^-)\otimes U^0\otimes U_\chi(\nnn).$$ 
The action of $\hhh$ stabilizes $J$ and induces an action on $U_\chi$. Set $(U_\chi)_{\mu+\bbz I}$ equal to the direct sum of weight spaces with weights in $\mu+\bbz I\in X/\bbz I$.
 We get therefore an $X/\bbz I$-grading on $U_\chi$. Note that there is a partial ordering relation on $X/\bbz I$ with $\lambda+\bbz I\leq \mu+\bbz I$ if and only if $\mu-\lambda\in\sum_{\alpha\in R^+\backslash \bbz I}\bbz^{\geq0}\alpha+\bbz I$. 
 The usual order relation on $X$ can also obtained from this relation if $I=\emptyset$. 

Let $A$ be a commutative Noetherian $k$-algebra over $U^0$ with the structural map $\pi\colon U^0\rightarrow A$.
We now define the category $\mathcal{C}_A(\ggg)$ whose objects are  $U_{\chi}\otimes A$-module $M$ with an $X/\bbz I$-grading satisfying the following conditions$\colon$
	\begin{enumerate}
		\item $M=\bigoplus_{\lambda+\bbz I\in X/\bbz I}M_{\lambda+\bbz I}$ as a $k$-vector space and each $M_{\lambda+\bbz I}$ are finitely generated over $A$.
		\item The action of $U_\chi$ and $A$ are commutative and compatible with the $X/\bbz I$-gradring, that is, $AM_{\lambda+\bbz I}\subset M_{\lambda+\bbz I}$ and $(U_{\chi})_{\mu+\bbz I}M_{\lambda+\bbz I}\subset M_{\lambda+\mu+\bbz I}$ for all $\lambda+\bbz I,\mu+\bbz I\in X/\bbz I$. 
		\item There is a weight decomposition $M_{\lambda+\bbz I}=\bigoplus_{d\mu\in\hhh^*,\mu\in\lambda+\bbz I+pX}M_{\lambda+\bbz I}^{d\mu}$ with $$h.m=(\pi(h)+d\mu(h))m$$
		for each $h\in\hhh$, and $m\in M_{\lambda+\bbz I}^{d\mu}$
	\end{enumerate}
	The morphsim of $\Ca(\ggg)$ are $U_\chi\otimes A$-homomorphisms perserving the $X/\bbz I$-grading. For convenience, we sometimes denote $\Ca$ instead of $\Ca(\ggg)$ if there is no confusion possible. This category is precisely the category in  \cite[3.2]{Ja00} and \cite[3.2]{We22}. It is compatible with the category in \cite[2.3/4]{AJS94} when $\chi=0$. 
	
	Suppose that $M\in\Ca$ and $m\in M_{\lambda+\bbz I}^{d\mu}$. Note that $(\pi(h_\alpha)+d\mu(h_\alpha))^p-(\pi(h_\alpha)+d\mu(h_\alpha))=\pi(h_\alpha)^p-\pi(h_\alpha)$. Then the central element $h_\alpha^p-h_\alpha$ acts as multiplication by $\pi(h_\alpha)^p-\pi(h_\alpha)$ on all $M_{\lambda+\bbz I}^{d\mu}$ and hence on $M$.
		
	When $A=F$ is a field, let $\xi$ be a linear function on $\ggg$ such that $\xi(h_\alpha)$ is in some purely inseparable extension of $F$ with $\xi(h_\alpha)^p=\pi(h_\alpha)^p-\pi(h_\alpha)$,
	 $\xi(x_{-\alpha})=1$ for $\alpha\in I$ and $\xi(x_{-\alpha})$=0 for $\alpha\in R\backslash{I}$.
	 Denote by $$U_\xi(\ggg\otimes F)=U(\ggg\otimes F)/(X^p-X^{[p]}-\xi(X)^p\mid X\in\ggg)$$
	 the reduced enveloping algebra of $\ggg\otimes F$ over $F$,
	 which is well-defined since $\xi(X)^p\in F$ by definition.
	 Then the category $\mathcal{C}_F$ is equivalent to the category of finite dimensional $X/\bbz I$-graded $U_\xi(\ggg\otimes F)$-modules 
     satisfy (2)(3) in definition.
	
	If $A=k$ and $\pi=0$, the category $\mathcal{C}
	_k$ is equivalent to a subcategory of $X/\bbz I$-graded $U_\chi(\ggg)$-module, cf.\cite[11.5]{Ja98}

    Let $M\in\Ca$. For $\nu\in X$, we define a grading shift object $M\langle\nu+\bbz I\rangle$ via defining $(M\langle\nu+\bbz I\rangle)_{\mu+\bbz I}=M_{\nu-\mu+\bbz I}$. 
    
    For any  $M\in\Ca$ and $\lambda\in X$, we have $x_\alpha.M_{\lambda+\bbz I}\subset M_{\lambda+\alpha+\bbz I}$. We get therefore each $$M^{[\lambda+\bbz R]}=\sum_{\mu\in\bbz R}M_{\lambda+\mu+\bbz I}$$
    is a $X/\bbz I$-graded submodule of $M$ and is an object in $\mathcal{C}_A$. We have $M=\bigoplus_{\nu\in X/\bbz R}M^{[\nu+\bbz R]}$. 
    Let $\Ca^{[\lambda+\bbz I]}$ be the full subcategory of $\Ca$ consisting of $M\in\Ca$ such that $M=M^{[\lambda+\bbz R]}$, and then
    we have $\Ca=\bigoplus_{\nu\in X/\bbz R}\Ca^{[\nu+\bbz I]}$.   
     Each $\Ca^{[\nu+\bbz I]}$ is equivalent to the category of finite dimensional $\bbz R/\bbz I$-graded $U_\chi\otimes A$-module such that $X^p-X^{[p]}$ acts as multiplication $\xi(X)^p$ and satisfy (2) in definition.
    
    Suppose $A=F$ is a field, and $\bbz R=\bbz I$, then
    each $\mathcal{C}_F^{[\nu+\bbz R]}$ 
     is equivalent to the category of finite dimensional $U_\xi(\ggg\otimes F)$-module. It follows that $\mathcal{C}_F$ are direct sum of some subcategories with each direct summand is equivalent to the category of finite dimensional $U_\xi(\ggg\otimes F)$-modules over $F$.

 Let $U_I$ be the image of $U_\chi(\ggg_I)$ in $U_\chi=U(\ggg)/J$ by regarding $U_\chi(\ggg_I)$ as a subalgebra of $U_\chi(\ggg)$. We also have an isomorphism of vector spaces by the PBW theorem
 $$U_I\simeq U_\chi(\nnn_I^-)\otimes U^0\otimes U_\chi(\nnn_I).$$

Since the $X/\bbz I$-grading on the subalgebra of $U_\chi$ can be obtained from the $X/\bbz I$-grading on $U_\chi$. Similarly, we can define several categories
\begin{enumerate}
	\item Define the category $\Ca(\hhh)$ to be the category of $U^0\otimes A$-module with analogous conditions.
	\item Define the category $\Ca(\bbb)$ to be the category of $U^0U_\chi(\nnn)\otimes A$-module with analogous conditions.
	\item Define the category $\Ca(\ggg_I)$ to be the category of $U_I\otimes A$-module with analogous conditions.
	\item Define the category $\Ca(\ppp)$ to be the category of $U_\chi(\nnn_I^-)U^0U_\chi(\nnn)\otimes A$-module with analogous conditions.
\end{enumerate}

Note that $(U_\chi)_{0+\bbz I}=U_I$. For any $M\in\Ca$ and $\lambda+\bbz I\in X/\bbz I$, each $M_{\lambda+\bbz I}$ is a $U_I\otimes A$-module and hence an object in $\mathcal{C}_{A}(\ggg_I)$. For any $f\in\Hom_{\Ca}(M,N)$, we have $f_{\lambda+\bbz I}\colon M_{\lambda+\bbz I}\rightarrow N_{\lambda+\bbz I}$ is a morphism in $\Ca(\ggg_I)$.

\begin{lem}\label{lem enough proj}
	Keep assumptions above.
	\begin{enumerate}
	\item There exists enough projectives in $\mathcal{C}_A(\ggg)$, $\mathcal{C}_A(\ggg_I)$, $\mathcal{C}_A(\bbb)$ and $\mathcal{C}_A(\ppp)$.
	\item Any projective object in $\mathcal{C}_A(\ggg)$, $\mathcal{C}_A(\ggg_I)$, $\mathcal{C}_A(\bbb)$ or $\mathcal{C}_A(\ppp)$ is still projective over $A$.
	\end{enumerate}
\end{lem}
 \begin{proof}
 	(1) followed by Lemma 4.2 and Corollary 4.13 in \cite{We22} and (2) followed by Lemma 2.7(c) in \cite{AJS94}.
 \end{proof}
\subsection{Subset $R_\pi$}
  We introduce a subset of $R$ which plays a key role in $\Ca$. 
  Let $A$ be a $U^0$-algebra with structural map $\pi\colon U^0\rightarrow A$. We denote
  $$R_\pi\colon=\{\alpha\in R\mid \pi(h_\alpha)^p-\pi(h_\alpha) \text{ is not a unit}\},$$ 
  and set $R_\pi^+=R^+\cap R_\pi$. We have $R_\pi=R^+_\pi\cup -R^+_\pi$.
   If $A=F$ is a field, then $R_\pi=\{\alpha\in R\mid \pi(h_\alpha)\in \mathbb{F}_p\}$.
	
	Set $B=U^0[(h_\alpha^{p-1}-1)^{-1}\mid \alpha\in R]$. For any $\alpha\in R^+$, denote $B^\alpha=B[h_\beta^{-1}\mid \beta\in R^+, \beta\neq\alpha]$ and $B^\emptyset=B[h_\beta\mid\beta\in R^+]$. We have the following observation, cf.\cite[5.3]{AJS94}.
\begin{lem}\label{lem R pi}
Let $A$ be a $U^0$-algebra with structural map $\pi\colon U^0\rightarrow A$.
\begin{enumerate}
\item $A$ is a $B^\emptyset$-algebra if and only if $R_\pi=\emptyset$.
\item For each $B$-algebra $A$ and $\alpha\in R^+$, we have $A$ is a $B^\alpha$-algebra if and only if $R^+_\pi\subset\{\alpha\}$.
\item If $A$ is a field. Then $A$ is a $B$-algebra if and only if $R_\pi=\{\alpha\in R\mid \pi(h_\alpha)=0\}$.
\end{enumerate}
\end{lem}
The assumption in \cite{We22} that $\pi(h_\alpha)=0$ for all $\alpha\in R_I$ is equivalent to $R_I\subset R_\pi$. When $A=k$ and $\pi=0$, we have $R_\pi=R$.
In the case of $A=F$ is a field. The subset $R_\pi\cap I$ determines the isomorphism class of irreducible object in $\mathcal{C}_F$ (see Corollary \ref{coro iso class of Z}) and $R_\pi\backslash R_I$ determines the irreducibility of some modules (see Corollary \ref{coro irreducible of Z}).

\section{Deformed baby Verma module}
In this section, we always assume that $A$ is a Noetherian commutative $U^0$-algebra with structural map $\pi\colon U^0\rightarrow A$. We keep assume $\chi\in\ggg^*$ has standard Levi form and that $I$ is the set of simple roots $\alpha$ with $\chi(x_\alpha)=0$ and $\xi$ is the corresponding linear form when $A=F$ is a field.
For any $A$-algebra $A'$, we shall assume it is Noetherian and commutative.

\subsection{Definition}
 For each $\lambda\in X$, let $A_\lambda\in\Ca(\hhh)$ denote the $U^0\otimes A$-module $A$ where each $h\in\hhh$ acts as multiplication by $d\lambda(h)+\pi(h)$. It can be extended into $U^0U_\chi(\nnn)\otimes A$-module via the projection $U^0U_\chi(\nnn)\twoheadrightarrow U^0$ and can be regarded as an object in $\Ca(\bbb)$. 
 Then we have induced module 
 $$Z_A(\lambda)=U_\chi\otimes_{U^0U_\chi(\nnn)}A_\lambda,$$
 which is called the deformed baby Verma module. (This name comes from \cite{So95}). 
 As a $U_\chi(\nnn^-)\otimes A$-module, we have an isomorphism $Z_A(\lambda)\simeq U_\chi(\nnn^-)\otimes A$.
 The $X/\bbz I$-grading on $Z_A(\lambda)$ is given by 
 $$Z_A(\lambda)_{\mu+\bbz I}=U_\chi(\nnn^-)_{\mu-\lambda+\bbz I}\otimes A_\lambda$$
 for any $\mu\in X$.
 We have $Z_A(\lambda)$ is an object in $\Ca$.

 Similarly, for each $\lambda\in X$, we can define the deformed baby Verma module in $\mathcal{C}_A(\ggg_I)$ as
 $$Z_{A,I}(\lambda)=U_I\otimes_{U^0U_\chi(\nnn_I)}A_\lambda,$$
 where $U_\chi(\nnn_I)$ acts as 0 on $A_\lambda$. 
 We can regard $Z_{A,I}(\lambda)$ as a object in $\Ca(\ppp)$ via the projection $\ppp\twoheadrightarrow \ggg_I$, and then we have
  $Z_A(\lambda)
  =U_\chi\otimes_{U_\chi(\nnn_I^-)U^0U_\chi(\nnn)}Z_{A,I}(\lambda)$. Since $(U_\chi)_{0+\bbz I}=U_I$, we also can get $Z_A(\lambda)_{\lambda+\bbz I}=Z_{A,I}(\lambda)$.
 
 When $A=F$ is a field. Consider as a $U_I\otimes F$-module, each deformed baby Verma module in $\mathcal{C}_F(\ggg_I)$ is isomorphic to some induced module 
 $$Z_{F,I}(d\lambda+\pi)=U_\xi(\ggg)\otimes_{U_\xi(\hhh\otimes\nnn_I)}F_\lambda,$$
 where $F_\lambda$ is 1-dimensional $U_\xi(\hhh\otimes\nnn_I)$-module with $h\in\hhh$ acts as multiplication by $d\lambda(h)+\pi(h)$ and $\nnn_I$ acts as zero. This module is extended from baby Verma module over an algebraically closed field. Since we have $\xi(h)^p=(d\lambda(h)+\pi(h))^p-(d\lambda(h)+\pi(h))=\pi(h)^p-\pi(h)$, the induced module is well-defined. 
 Note that $\chi|_{\ggg_I}$ is regular nilpotent, we have $U_\chi(\nnn_I^-)\otimes F$ is simple over itself. Since we have an isomorphism of $U_\chi(\nnn_I^-)\otimes F$-modules $Z_{F,I}(\lambda)\simeq U_\chi(\nnn_I^-)\otimes F$, 
  each $Z_{F,I}(\lambda)$ is irreducible for  $\lambda\in X$.

\begin{lem}\label{Z maximal submodule}
Let $A=F$ be a field, then $Z_F(\lambda)$ has a unique maximal submodule.
\end{lem}
\begin{proof}
Note that $(U_\chi(\nnn^-)\otimes F)_{<0+\bbz I}=\bigoplus_{\nu+\bbz I<0+\bbz I}(U_\chi(\nnn^-)\otimes F)_{\nu+\bbz I}$ is an ideal in $U_\chi(\nnn^-)\otimes F$ with quotient
$(U_\chi(\nnn^-)\otimes F)_{0+\bbz I}=U_\chi(\nnn_I^-)\otimes F$. 
Since $U_\chi(\nnn_I^-)\otimes F$ is a simple algebra, then $(U_\chi(\nnn^-)\otimes F)_{<0+\bbz I}$ is the radical of $U_\chi(\nnn^-)\otimes F$. 
Under the isomorphism $Z_F(\lambda)\simeq U_\chi(\nnn^-)\otimes F$ as a $U_\chi(\nnn^-)\otimes F$-module, every proper submodules of $Z_F(\lambda)$ are contained in the subspace $(U_\chi(\nnn^-)\otimes F)_{<0+\bbz I}\otimes F$. Hence $Z_F(\lambda)$ has a unique maximal submodule.
\end{proof}
When $A=F$ is a field, we denote the unique irreducible quotient of $Z_F(\lambda)$ as 
$L_F(\lambda)=Z_F(\lambda)/\Rad Z_F(\lambda)$.
  We have obviously $\Rad Z_F(\lambda)$ and $L_F(\lambda)$ are both in $\mathcal{C}_F$. Each simple module in $\mathcal{C}_F$ is of form $L_F(\lambda)$ with $\lambda\in X$.   
Note that the surjective map $Z_F(\lambda)\twoheadrightarrow L_F(\lambda)$ induce a nonzero map $Z_F(\lambda)_{\lambda+\bbz I}=Z_{F,I}(\lambda)\rightarrow L_F(\lambda)_{\lambda+\bbz I}$ in $\mathcal{C}_F(\ggg_I)$.
The simplicity of $Z_{F,I}(\lambda)$ shows that  $Z_{F,I}(\lambda)=L_F(\lambda)_{\lambda+\bbz I}$.

In the case of $A$ being a ring and $\mmm$ is a maximal ideal in $A$, the irreducible object $L_{A/\mmm}(\lambda)$ can be regarded as an object in $\Ca$ which is annihilated by $\mmm$. The irreducible object is indexed by $X$ and the maximal ideals in $A$, see similar results in \cite[1.7.3]{GJ81}.

Different from the restricted case, not all $Z_F(\lambda)$ are different objects in $\mathcal{C}_F$ for $A=F$ is a field. We have to answer when two deformed baby Verma modules are isomorphic. 

We first see the coadjoint action of $G$ in linear form $\xi$. Suppose $A=F$ is a field. Since $(\xi(x)+\xi(y))^p=\xi(x)^p+\xi(y)^p$ for any $x,y\in\ggg$. We have $(g\cdot\xi(x))^p\in F$ for any $g\in G$ and $x\in\ggg$. It follows that $U_{g\cdot\xi}(\ggg\otimes F)$ is well defined and we have $U_\xi(\ggg\otimes F)\simeq U_{g\cdot\xi}(\ggg\otimes F)$.

A decomposition $\xi=\xi_s+\xi_n$ is called Jordan decomposition, if there exists $g\in G$ such that $g\cdot\xi_s(\nnn\oplus\nnn^-)=0$, $g\cdot\xi_n(\bbb)=0$ and $g\cdot\xi_n(x_\alpha)=0$ if $g\cdot\xi_s(h_\alpha)\neq0$. 
Let $\ggg_\xi=\{x\in\ggg\mid \xi([X,\ggg]=0)\}$ be the stabilizer of $\xi$ in $\ggg$. 
We call $\xi$ is regular if there exists $g\in G$ such that $g\cdot\xi$ is regular nilpotent over $[\ggg_{(g\cdot\xi)_s},\ggg_{(g\cdot\xi)_s}]$, cf. \cite[4]{FP88}. It is obvious that regular nilpotent linear form is regular.
%
 
Denote  $W_{I,\pi}=\langle s_{\alpha}\in W_I\mid \alpha\in R_\pi\cap I\rangle$ be the subgroup of $W_I$.
\begin{lem}\label{propertty of g_I}
Suppose $B$-algebra $F$ is a field. Let $R^+_\pi\cap R_I$ equal to $\emptyset$, $\{\alpha\}$, or $R^+_I$. Then
\begin{enumerate}
\item $Z_{F,I}(\lambda)\simeq Z_{F,I}(\mu)$ if and only if $\lambda\in W_{I,\pi}\cdot\mu+p\bbz I$;
\item $\Ext_{\mathcal{C}_F}^{\bullet}(Z_{F,I}(\lambda),Z_{F,I}(\mu))\neq0$ if and only if $\lambda\in W_{I,\pi}\cdot\mu+p\bbz I$;
\item The projective cover $Q_{F,I}(\lambda)$ of $Z_{F,I}(\lambda)$ is of length $|W_{I,\pi}\cdot \lambda|$ with all composition factors ismorphic to $Z_{F,I}(\lambda)$.
\end{enumerate}
\end{lem}
\begin{proof}
Since $\mathcal{C}_F(\ggg_I)$ is the direct sum of full subcategories with each isomorphic to the category of finite dimensional $U_\xi(\ggg_I\otimes F)$-modules and $Z_{F, I}(\lambda)\simeq Z_{F, I}(d\lambda+\pi)$ as $U_\xi(\ggg_I)$-modules. 
So $Z_{F,I}(\lambda)\simeq Z_{F,I}(\mu)$ if and only if $Z_{F,I}(d\lambda+\pi)\simeq Z_{F,I}(d\mu+\pi)$ and $\lambda+\bbz I=\mu+\bbz I$.

In case of $R^+_\pi\cap R_I=\emptyset$ or $R^+_\pi \cap R_I=\{\alpha\}$ with $\alpha\notin I$, we have $(d\lambda+\pi)(h_\beta)\notin \mathbb{F}_p$ for any simple root $\beta$, and $Z_{F,I}(d\lambda+\pi)$ is projective $U_\xi(\ggg_I\otimes F)$-module, cf.\cite[Theorem 4.2]{FP90}. Note that $W_{I,\pi}=1$, then (3) follows. Note that 
$$\dim Z_{F,I}(d\lambda+\pi)=p^{\dim\nnn_I}\dim Z_{F,I}(d\lambda+\pi)^{\nnn_I},$$ cf.\cite[10.9]{Ja98}. 
We have $\dim Z_{F,I}(d\lambda+\pi)^{\nnn_I}=1$, and then $Z_{F,I}(d\lambda+\pi)\simeq Z_{F,I}(d\mu+\pi)$ if and only if $d\lambda=d\mu$. Since $pX\cap\bbz I=p\bbz I$, we have $\lambda\in\mu+p\bbz I$ and (1)(2) follows.

For the remaining cases, by comparing \cite[Theorem 4.4]{FP88}, it is enough to show that $\xi$ is regular. 
If $R^+_\pi\cap R_I=R^+_I$, we have $\xi(h_\alpha)=0$ for all $\alpha\in R_I$ and hence $\xi$ is a regular nilpotent linear form on $\ggg_I$. 
If $R^+_\pi\cap I=\{\alpha\}$ with $\alpha\in I$, then for any $\beta\in R_I^+$, $\pi(h_\beta)=0$ only when $\beta=\alpha$.  
Note that $-\alpha+\beta$ is not a negative root for any $\beta\in R_I^+\backslash R^+_\pi$. 
By computation in \cite[3.8]{KW76},
there exists $g\in G$ such that $(g\cdot\xi)(\oplus_{\beta\neq-\alpha}\ggg_{\beta})=0$, $(g\cdot\xi)(x_{-\alpha}\neq0)$, $(g\cdot\xi)(h_\alpha)=0$, and $(g\cdot\xi)(h_\beta)\neq0$ for any $\beta\neq\pm\alpha$. 
It follows that $g\cdot\xi$ is regular nilpotent on $[\ggg_{(g\cdot\xi)_s},\ggg_{(g\cdot\xi)_s}]=\ggg_{-\alpha}\oplus\hhh\oplus\ggg_{\alpha}$. 
Then the linear form $\xi$ is regular as desired.
\end{proof}

\begin{coro}\label{coro iso class of Z}
Suppose $B$-algebra $F$ is a field and let $\lambda,\mu\in X$. Then the following are equivalent$\colon$ 
\begin{enumerate}
	\item $Z_F(\lambda)\simeq Z_F(\mu)$.
	\item $L_F(\lambda)\simeq L_F(\mu)$.
	\item $\lambda\in W_{I,\pi}\cdot\mu+p\bbz I$.
\end{enumerate}
\end{coro}
\begin{proof}
If $Z_F(\lambda)\simeq Z_F(\mu)$, then  $L_F(\lambda)\simeq L_F(\mu)$ is obviously. Suppose $L_F(\lambda)\simeq L_F(\mu)$, then we have $\lambda+\bbz I=\mu+\bbz I$ as the largest nonzero $X/\bbz I$-grading of homogeneous subspace. It follows that
 $L_F(\lambda)_{\lambda+\bbz I}\simeq L_F(\mu)_{\lambda+\bbz I}$, and then we have $Z_{F,I}(\lambda)\simeq Z_{F,I}(\mu)$. We get therefore $\lambda\in W_{I,\pi}\cdot\mu+p\bbz I$ by Lemma \ref{propertty of g_I}. 
Finally, 
suppose $\lambda\in W_{I,\pi}\cdot\mu+p\bbz I$. We have $Z_{F}(\lambda)\simeq Z_F(\mu)$ by \cite[Proposition 6.2]{We22}.
\end{proof}
\subsection{Properties}
We recall and generalize some results in \cite{AJS94} and \cite{We22}. It will be useful in the later sections.
\begin{lem}\label{lem of Ext}
	Let $\lambda, \mu\in X$.
	\begin{enumerate}
		\item If $\Ext^1_{\Ca}(Z_A(\lambda),M)\neq 0$ for some $M\in\Ca$, then there exists $\mu+\bbz I\in X/\bbz I$ such that $M_{\mu+\bbz I}\neq0$ and $\lambda+\bbz I\leq\mu+\bbz I$. In particular, if $\lambda+\bbz I\nleq\mu+\bbz I$, then $\Ext^1_{\mathcal{C}_A}(Z_A(\lambda),Z_A(\mu))=0$.
		\item If $\lambda+\bbz I=\mu+\bbz I$, then $\Ext^1_{\mathcal{C}_A}(Z_A(\lambda),Z_A(\mu))=0$
		when $\Ext^1(Z_{A,I}(\lambda),Z_{A,I}(\mu))=0$.	
		\item Suppose $A=F$ is a field. If $\lambda+\bbz I\nleq\mu+\bbz I$, then
		$$\Ext^1_{\mathcal{C}_F}(L_F(\lambda),L_F(\mu))\simeq\Hom_{\mathcal{C}_F}(\Rad Z_F(\lambda),L_F(\mu)).$$
	\end{enumerate}
\end{lem}	
\begin{proof}
	(1)See \cite[5.2]{We22}.
	
	(2)Suppose that $\lambda+\bbz I=\mu+\bbz I$. Consider an arbitrary extension of $Z_A(\lambda)$ and $Z_A(\mu)$ in $\Ca\colon$
	$$0\rightarrow Z_A(\lambda)\longrightarrow M\longrightarrow Z_A(\mu)\rightarrow0.$$	
	Note that $Z_A(\lambda)_{\lambda+\bbz I}=Z_{A,I}(\lambda)$ and $Z_A(\mu)_{\lambda+\bbz I}=Z_A(\mu)_{\mu+\bbz I}=Z_{A,I}(\mu)$. If we restrict the extension to the subspace of degree of $\lambda+\bbz I$, it yields an extension
	$$0\rightarrow Z_{A,I}(\lambda)\rightarrow M_{\lambda+\bbz I}\rightarrow Z_{A,I}(\mu)\rightarrow 0$$ 
	in $\Ca(\ggg_I)$. Note that the generators of $Z_A(\lambda)$ and $Z_A(\mu)$ is contained in $Z_A(\lambda)_{\lambda+\bbz I}$ and $Z_A(\mu)_{\lambda+\bbz I}$. It follows that $\Ext^1_{\mathcal{C}_A}(Z_A(\lambda),Z_A(\mu))=0$ when $\Ext^1(Z_{A,I}(\lambda),Z_{A,I}(\mu))=0$.
	
	(3) Suppose $\lambda+\bbz I\nleq \mu+\bbz I$ for $\lambda, \mu\in X$. We have a short exact sequence 
	$$0\rightarrow\Rad Z_F(\lambda)\longrightarrow Z_F(\lambda)\longrightarrow L_F(\lambda)\rightarrow0.$$
	By applying the functor $\Hom_{\mathcal{C}_F}(-,L_F(\mu))$ for any $\mu\in X$ to this exact sequence and the fact $Z_F(\lambda)$ has a simple head $L_F(\lambda)$,
	we get an exact sequence
	$$0\rightarrow\Hom(\Rad Z_F(\lambda),L_F(\mu))\longrightarrow\Ext^1(L_F(\lambda),L_F(\mu))\longrightarrow\Ext^1(Z_F(\lambda),L_F(\mu)).$$
	Since $\lambda+\bbz I\nleq\mu+\bbz I$, the last term is zero by (1) and the isomorphism follows.	
\end{proof}

A $Z$-filtration of $M\in\Ca$ is a chain of submodule
 $$M_0=0\subset M_1\subset M_2\cdots\subset M_r=M,$$
 such that each $M_{i}/M_{i-1}$ isomorphic to some $Z_A(\lambda_i)$ with $\lambda_i\in X$. 

\begin{lem}\label{lem proj Z-fil}Let $M\in\Ca$.
\begin{enumerate}
\item There exists a projective object $Q$ in $\Ca$ with a $Z$-filtration such that $M$ is a homomorphic image of $Q$. 
\item Suppose $A$ is a local ring, then any projective object in $\Ca$ has a $Z$-filtration.
\end{enumerate}
\end{lem}
\begin{proof}
(1)See \cite[Theorem 4.10]{We22}.

(2)Let $M$ be a projective module in $\Ca$. Then there exists a projective object $Q$ with an epimorphism from $Q$ onto $M$. 
By the projectivity of $M$, the epimorphism has to split and we get therefore $M$ is a direct summand of $Q$.
Since we have a standard argument over local ring$\colon$ Each direct summand of a module with a $Z$-filtration has a $Z$-filtraion, cf.\cite[2.16]{AJS94},
then $M$ has a $Z$-filtration as desired. 
\end{proof}

Suppose $A'$ is an $A$-algebra, then $A'$ is naturally a $U^0$-algebra. From \cite[3.1]{AJS94} or \cite[5.3]{We22}, we have extension of scalars functor $-\otimes_A A'\colon\Ca\rightarrow\mathcal{C}_{A'}$ which sends $M$ to $M_{A'}=M\otimes_A A'$. We have the following results.

\begin{lem}\label{lem base change}
Let $M\in\Ca$ and $A'$ is a $A$-algebra.
\begin{enumerate}
\item If $M$ is projective, so is $M\otimes_A A'$ in $\mathcal{C}_{A'}$.
\item If $M$ has a $Z$-filtration, then $M\otimes_A A'$ also has a $Z$-filtration.
\item If $A'$ is flat over $A$. We have 
$$\Ext^i_{\Ca}(M,N)\simeq \Ext^i_{\mathcal{C}_{A'}}(M\otimes_A A',N\otimes_A A')$$
for all $M,N$ in $\Ca$ and $i\geq0$.
\end{enumerate}	
\end{lem}
\begin{proof}
	(1) and (2) followed by Lemma 5.7 and Lemma 5.8 in \cite{We22}. (3) comes from \cite[Lemma 3.2]{AJS94}.
\end{proof}
Suppose $A=F$ is a field. 
Since $\mathcal{C}_F$ has enough projective object, we debote by $Q_F(\lambda)$ the unique projective cover of $L_F(\lambda)$ in $\mathcal{C}_F$ for each $\lambda\in X$.

\begin{lem}\label{lem proj over local ring}
	Suppose $A$ is a local ring with resuid field $F$ and $R_I\subset R_\pi$. Then there exists a projective object $Q_A(\lambda)$ in $\Ca$ such that $Q_A(\lambda)\otimes F=Q_{F}(\lambda)$.
\end{lem}
\begin{proof}
	Since $R_I\subset R_\pi$, it yields that $\pi(h_\alpha)=0$ for any $\alpha\in R_I$. Then results follows from \cite[Theorem 8.29]{We22}.
\end{proof}
We have the following results, which help us focus on the case of the field.
\begin{lem} \label{lem ext base change}
	Let $M$, $N$ be modules in $\Ca$ that are projective over $A$. Let $A'$ be an $A$-algebra with finite projection resolution as an $A$-module. Suppose
	$$\Ext^i_{\mathcal{C}_{A/\mmm}}(M\otimes A/\mmm, N\otimes A/\mmm)=0$$ for all maxiaml ideals $\mmm$ of $A$ and all $i>0$. Then 
	\begin{enumerate}
		\item $\Hom_{\Ca}(M,N)$ is a projective $A$-module.
		\item $\Ext^i_{\Ca}(M\otimes A', N\otimes A')=0$ for all $i>0$
		\item $\Hom_{\Ca}(M,N)\otimes A'=\Hom_{\Ca}(M\otimes A', N\otimes A')$
	\end{enumerate}
\end{lem}   
\begin{proof}
	Since $A$ is Noetherian, then $A_\mmm$ is a regular local ring. The results follow from \cite[3.4]{AJS94}.
\end{proof}  
\begin{coro}\label{coro projective from field to ring}
	Suppose $M$ is a module in $\Ca$. Then $M$ is a projective object in $\Ca$ if and only if $M\otimes_A A/\mmm$ is projective in $\mathcal{C}_{A/\mmm}$ for all maxiaml ideals  $\mmm$ of $A$.
\end{coro}  
\begin{proof}
Similar to \cite[Corollary 3.5]{AJS94}.
\end{proof}
\subsection{Twists and homomorphisms}\label{sec twist}
For any $w\in W$, let $w\nnn=\oplus_{\alpha\in R^+}\ggg_{w\alpha}$ and set $w\bbb=\hhh\oplus w\nnn$. We have $\chi(w\bbb)=0$ if and only if $w\in W^I$, where $W^I=\{w\in W\mid w^{-1}\alpha>0 \text{ for all } \alpha\in I\}$. The PBW theorem shows that we have an isomorphism of vector spaces
$$U_\chi\simeq U_\chi(w\nnn^-)U^0U_\chi(w\nnn).$$
Similarly, for $w\in W^I$, we can define $U^0U_\chi(w\nnn)\otimes A$-module $A_\lambda$ where $U_\chi(w\nnn)$ acts as 0 and each $h\in\hhh^*$ act as multiplication $d\lambda(h)+\pi(h)$. We can define  twist deformed baby Verma module
$$Z_A^{w}(\lambda)=U_\chi\otimes_{U^0U_\chi(w\nnn)}A_\lambda,$$
for each $\lambda\in X$ and $w\in W^I$. 
There is a unique $X/\bbz I$-grading on $Z_A^{w}(\lambda)$ such that $Z_A^{w}(\lambda)_{\mu+\bbz I}=U_\chi(w\nnn^-)_{\mu-\lambda+\bbz I}\otimes A_\lambda$ and $Z^w_A(\lambda)$ is an object in $\Ca$. When $w=\id$, we have $Z^{\id}_A(\lambda)=Z_A(\lambda)$. 
Since $w\nnn_I^-=\nnn_I^-$ for any $w\in W^I$, we have $Z_A^{w}(\lambda)_{\lambda+\bbz I}=Z_{A,I}(\lambda)$ in $\Ca(\ggg_I)$.
If $A=F$ is a field, we can generalize Lemma \ref{Z maximal submodule} to $Z^{w}_F(\lambda)$. So we have the irreducible quotient $L^{w}_F(\lambda)=Z^{w}_F(\lambda)/\Rad Z^w_A(\lambda)$ in $\mathcal{C}_F$ for each $w\in W^I$.

We have following lemma from \cite{Ja00}.
\begin{lem}\label{lemma of homo}
Suppose $\lambda\in X$. Let $w\in W^I$ and $\alpha$ being a simple root such that $w\alpha>0$ and $ws_\alpha\in W^I$. In the category $\Ca$, there exists homomorphsims $\varphi\colon Z^w_{A}(\lambda-(p-1)w\alpha)\longrightarrow Z^{ws_\alpha}_{A}(\lambda)$, given by $\varphi(1\otimes 1)=x_{w\alpha}^{p-1}\otimes 1$ and 
	$\varphi'\colon Z^{ws_\alpha}_{A}(\lambda)\longrightarrow Z^w_{A}(\lambda-(p-1)w\alpha)$, given by $\varphi'(1\otimes 1)=x_{-w\alpha}^{p-1}\otimes 1$.
\end{lem}

For $M\in\Ca$, we write $[M]$ for its image in the Grothendieck group of category $\Ca$. Note that $w_Iw_0$ is the unique longest element in $W^I$ and we denote it by $w^I$.
\begin{prop}\label{twist homo property}
Keep the assumptions and notations in Lemma \ref{lemma of homo}. Suppose $A=F$ is a field.
	\begin{enumerate}
	\item Suppose $w\alpha\notin R_\pi$, then $\varphi$ and $\varphi'$ are isomorphisms.
	\item Suppose $w\alpha\in R_\pi$. 
	Let $\pi(h_{w\alpha})+\langle\lambda,w\alpha^\vee\rangle+1\equiv d$ $(\mod p)$ with $0\leq d<p$. 
	Then we have 
 $$[Z^w_F(\lambda)]=[Z^{ws_\alpha}_F(\lambda-(p-1)w\alpha)].$$
     Furthermore, we have
	 $\varphi$ and $\varphi'$ are isomorphisms when $d=0$. And there exists a homomorphsim $\psi\colon Z^w_F(\lambda-dw\alpha)\rightarrow Z_F^w(\lambda)$ such that $\im\psi\simeq\ker\varphi$ when $d\neq0$.
	\end{enumerate}
\end{prop}
\begin{proof}
	Let $\beta=w\alpha$ and set
	 $$a_i=\dfrac{x^i_{-\beta}}{i!}\otimes 1\in Z_F^w(\lambda) \text{ and } b_i=\dfrac{x^i_{\beta}}{i!}\otimes 1\in Z^{ws_\alpha}_F(\lambda)$$
	 for $0\leq i<p$.
	Let $\uuu^-=\bigoplus_{\gamma\in R^+\backslash\{\alpha\}}\ggg_{-\gamma}$. 
	We have $U_\chi(w\nnn^-)\otimes F=\bigoplus_{i=0}^{p-1} U_\chi(w\uuu)a_i\otimes F$ and $U_\chi(ws_\alpha\nnn^-)\otimes F=\bigoplus_{i=0}^{p-1} U_\chi(w\uuu)b_i\otimes F$. 
	To determine homomorphisms $\varphi$ and $\varphi'$, we only need to check all  $\varphi(a_i)$ and $\varphi'(b_i)$.
	
	Note that
	$x_\beta a_i=(\pi(h_\beta)+d\lambda(h_\beta)-i+1)a_{i-1}$ and $x_{-\beta}b_i=-(\pi(h_\beta)+d\lambda(h_\beta)+i+1)b_{i-1}$. 
	We get therefore
	$$\varphi(a_i)=(-1)^i\dfrac{(p-1)!}{i!}\prod_{j=1}^{i}(\pi(h_\beta)+d\lambda(h_\beta)-j+1)b_{p-1-i},$$
	and
	$$\varphi(b_i)=\dfrac{(p-1)!}{i!}\prod_{j=1}^{i}(\pi(h_\beta)+d\lambda(h_\beta)+j+1)a_{p-1-i}.$$
		
	When $\beta\notin R_\pi$, we have $\pi(h_\beta)^p-\pi(h_\beta)$ is a unit, and then $\pi(h_{\beta})\notin\bbf_p$. Then $\pi(h_\beta)+\lambda(h_\beta)+k$ is nonzero for any $k\in\bbz$ and hence $\varphi$ and $\varphi'$ are isomorphisms by comparing elements. We get (1) as desired.
	
	When $\beta\in R_\pi$ with $d=0$. Then $\pi(h_\beta)+d\lambda(h_\beta)\pm j+1\neq0$ for all $0\leq j< p$ and hence $\varphi$ and $\varphi'$ are isomorphisms.
	
	When $\beta\in R_\pi$ with $d\neq0$. We have $\varphi(a_i)=0$ if and only if $i\geq d$. We get therefore
	$$\ker\varphi=\bigoplus_{i=d}^{p-1}(U_\chi(w\uuu^-)\otimes F)a_i\text{ and } \im\varphi=\bigoplus_{i=p-d}^{p-1}(U_\chi(w\uuu^-)\otimes F)b_i.$$
	Also, we have $\varphi'(b_i)=0$ if and only if $i\geq p-d$. We get therefore
	$$\ker\varphi'=\bigoplus_{i=p-d}^{p-1}(U_\chi(w\uuu^-)\otimes F)b_i \text{ and } \im\varphi'=\bigoplus_{i=d}^{p-1}(U_\chi(w\uuu^-)\otimes F)a_i.$$
	Hence, we get therefore $$[Z^w_F(\lambda)]=[\ker\varphi]+[\im\varphi]=[\im\varphi']+[\ker\varphi']=[Z^{ws_\alpha}_F(\lambda-(p-1)\beta]$$
 as desired. Note that $\ker\varphi\simeq ((U_\chi(w\nnn^-)\otimes F).a_d$ and $((U_\chi(w\nnn)\otimes F).a_d=0$. Then we have a homomorphsim $\psi\colon Z_F^{w}(\lambda-d\beta)\rightarrow Z_F^{w}(\lambda)$ with $\psi(1\otimes 1)=a_d$ such that $\im\psi\simeq\ker\varphi$.
\end{proof}
 
The deformed baby Verma module has a simple head. It also has a simple socle and can be described by twist deformed baby Verma module as follows. We put $\lambda^w=\lambda-(p-1)(\rho-w\rho)$ for any $w\in W^I$.
\begin{lem}\label{Z socle}
Suppose that $A=F$ is a field.
\begin{enumerate} 	
		\item $U_\chi(\nnn^-)\otimes F$ as a module over itself has a socle equal to $(U_\chi(\nnn^-)\otimes F)_{-(p-1)(\rho-w^I\rho)+\bbz I}$.
		\item  Let $\lambda\in X$, the scole of $Z_F(\lambda)$ or $Z^{w^I}_F(\lambda)$ is simple and equal to 
		$$\Soc(Z_F(\lambda))=U_\chi Z_F(\lambda)_{\lambda^{w^I}+\bbz I},$$ 
		and
		$$\Soc(Z^{w^I}_F(\lambda^{w^I}))=U_\chi  Z^{w^I}_F(\lambda^{w^I})_{\lambda+\bbz I}.$$	
	\end{enumerate}	
\end{lem}
\begin{proof}
	(1)
	Note that $(U_\chi(\nnn^-)\otimes F)_{-(p-1)(\rho-w^I\rho)+\bbz I}$ is annihiated by $(U_\chi(\nnn^-)\otimes F)_{<0+\bbz I}$, and is stable under the action of $(U_\chi(\nnn^-)\otimes F)_{0+\bbz I}$, so $(U_\chi(\nnn^-)\otimes F)_{-(p-1)(\rho-w^I\rho)+\bbz I}$ is a module over $U_\chi(\nnn^-)\otimes F$. 
	Since $(U_\chi(\nnn^-)\otimes F)_{0+\bbz I}=U_\chi(\nnn_I^-)\otimes F$ is simple over itself, then $(U_\chi(\nnn^-)\otimes F)_{-(p-1)(\rho-w^I\rho)+\bbz I}$ is simple $U_\chi(\nnn^-)\otimes F$-module and contained in every submodule of $U_\chi(\nnn^-)\otimes F$.
	Hence it equal to socle.
	
	(2)Since $Z_F(\lambda)$ isomorphic to $U_\chi(\nnn^-)\otimes F$ as a $U_\chi(\nnn^-)\otimes F$-module. We deduce from (1) that $Z_F(\lambda)_{\lambda^{w^I}+\bbz I}$ contained in every nonzero $U_\chi\otimes F$-submodule of $Z_F(\lambda)$. Then $Z_F(\lambda)$ has a simple socle, which lies in $\mathcal{C}_F$ obviously. It follows that $$\Soc(Z_F(\lambda))=U_\chi((U_\chi(\nnn^-)_{-(p-1)(\rho-w^I\rho)+\bbz I}\otimes F)=U_\chi Z_F(\lambda)_{\lambda^{w^I}+\bbz I}.$$
	On the other hand, note that $w^I(\nnn^-)=\nnn_I^-\bigoplus\bigoplus_{\beta\in R^+\backslash R^+_I}\ggg_\beta$. 
	We can similarly get $U_\chi(w^I\nnn^-)\otimes F$ as a module over itself has a socle equal to $(U_\chi(w^I\nnn^-)\otimes F)_{(p-1)(\rho-w^I\rho)+\bbz I}$. And then we have
	$$\Soc(Z^{w^I}_F(\lambda^{w^I}))=U_\chi  Z^{w^I}_F(\lambda^{w^I})_{\lambda+\bbz I}.$$

\end{proof}
\begin{prop}\label{prop long homo}
 Let $\lambda\in X$. We have homomorphisms
 $$\varphi\colon Z_A(\lambda)\rightarrow Z_A^{w^I}(\lambda^{w^I}),$$
 and
 $$\varphi'\colon Z_A^{w^I}(\lambda^{w^I}) \rightarrow Z_A(\lambda).$$
 Furthermore, if $A=F$ is a field, then $\im\varphi=L_F(\lambda)$ and $\im\varphi'=L_F^{w^I}(\lambda^{w^I})$.
\end{prop}
\begin{proof}
	Suppose $w^I=s_1s_2\cdots s_m$ is a reduction expression, where $s_i=s_{\alpha_i}$ with $\alpha_i$ is a simple root. Set $w_i=s_1s_2\cdots s_{i-1}$ for $1\leq i\leq m+1$. So $w_1=1$ and $w_{m+1}=w^I$. Note that
	 $$\lambda^{w_{i+1}}=\lambda-(p-1)(\rho-w_{i+1}\rho)=\lambda-(p-1(\rho-w_i\rho-\alpha))=\lambda^{w_i}-(p-1)\alpha,$$ 
	we have for $1\leq i\leq m+1$ homomorphisms $\varphi_i\colon Z^{w_i}_A(\lambda^{w_i})\rightarrow Z^{w_{i+1}}_A(\lambda^{w_{i+1}})$ and  $\varphi_i'\colon Z^{w_{i+1}}_A(\lambda^{w_{i+1}})\rightarrow Z^{w_i}_A(\lambda^{w_i})$.
	We get therefore the composition homomorphism
	$$\varphi=\varphi_m\circ\cdots\circ\varphi_1\colon Z_A(\lambda)=Z_A^{w_1}(\lambda^{w_1})\rightarrow Z_A^{w_2}(\lambda^{w_2})\rightarrow\cdots\rightarrow Z_A^{w_{m+1}}(\lambda^{w_{m+1}})=Z_A^{w^I}(\lambda^{w^I}),$$
	and
		$$\varphi'=\varphi'_1\circ\cdots\circ\varphi'_m\colon Z_A^{w^I}(\lambda^{w^I})=Z_A^{w_{m+1}}(\lambda^{w_{m+1}})
		\rightarrow Z_A^{w_m}(\lambda^{w_m})\rightarrow\cdots\rightarrow Z^{w_1}_A(\lambda^{w_1})=Z_A(\lambda).$$
Suppose now $A=F$ is a field. Since $Z_F(\lambda)$ is generated by $Z_F(\lambda)_{\lambda+\bbz I}$, the image of $\varphi$ is generated by $$\varphi(Z_F(\lambda)_{\lambda+\bbz I})=Z^{w^I}_F(\lambda^{w^I})_{\lambda+\bbz I}.$$ 
	But the right side is the generator of the simple socle by Lemma \ref{Z socle}. Hence, we have $\im\varphi$ is irreducible. Since $\im\varphi$ is isomorphic to the quotient of $Z_F(\lambda)$, we conclude that $\im\varphi\simeq L_F(\lambda)$.
	
	Similarly, we can get $\im\varphi'\simeq L_F^{w^I}(\lambda^{w^I})$.
\end{proof}

When $A=F$ is a field, for each $\beta\in R_\pi$, we denote by $n_\beta(\lambda)$ the integer such that $\langle\lambda+\rho,\beta^\vee\rangle+\pi(h_\beta)\equiv n_\beta(\lambda)\;(\mod p)$ with $0\leq n_\beta(\lambda)\leq p-1$.
\begin{coro} \label{coro irreducible of Z}Suppose $A=F$ is a field.
If $n_\beta(\lambda)=0$ for all $\beta\in R_\pi\backslash R_I$, then $Z_F(\lambda)$ is irreducible for all $\lambda\in X$. 
\end{coro}
\begin{proof}
When $\bbz I=\bbz R$, then $I$ is equal to the set of all simple roots. Then $Z_F(\lambda)=Z_{F,I}(\lambda)$ is simple for each $\lambda\in X$. 

Suppose $\bbz I\neq\bbz R$. We have a nontrivial homomorphism  $Z_F(\lambda)\rightarrow Z^{w^I}_F(\lambda^{w^I})$ from Proposition \ref{prop long homo}.
For any $w\in W^I$ and simple root $\alpha$, we have
$$\langle\lambda+\rho,w\alpha^\vee\rangle=\langle\lambda^{w},w\alpha^\vee\rangle+\langle\rho+(p-1)(\rho-w\rho),w\alpha^\vee\rangle\equiv\langle\lambda^w,w\alpha^\vee\rangle+1 \;(\mod p).$$
 By Proposition \ref{twist homo property}, the composition homomorphism $\varphi\colon Z_F(\lambda)\rightarrow Z^{w^I}_F(\lambda^{w^I})$ is isomorphism if and only if all $n_\beta(\lambda)=0$ for all $\beta\in R_\pi\backslash R_I$.
Hence, we have $Z_F(\lambda)\simeq\im\varphi\simeq L_F(\lambda)$ as desired.
\end{proof}
When $I=\emptyset$, it is just the result in \cite[6.3]{AJS94}. The subset $R_\pi\backslash R_I$ has deep connections with the irreducibility of the deformed baby Verma module. However, it is not projective and injective in general like that in the case of $I=\emptyset$, since two deformed baby Verma modules in $\Ca$ with different highest weights in $X$ may be isomorphic. We will discuss this in Section \ref{subsec proj} later.
\subsection{Duality}
We recall some constructions in \cite[5.4]{We22} and generalize some results in \cite[11.14/15]{Ja98}.

Let $A$ be a $U^0$-algebra with structural map $\pi\colon U^0\rightarrow A$. Denote the category of $U_{-\chi}\otimes A$-module with analogue conditions of $\Ca$ as $\overline{\Ca}$. We could similarly define the deformed baby Verma modules $\overline{Z}_A(\lambda)$ and their twist $\overline{Z}^{w}_A(\lambda)$ for $\lambda\in X$ and $w\in W^I$.
Write $\overline{A}$ the same $k$-algebra structure with $A$ with structural map $\overline{\pi}\colon U^0 \rightarrow {A}$ such that $\overline{\pi}(h)=-\pi(h)$ for $h\in U^0$.

Consider the contravariant dual of $M$ in $\Ca$ by setting
$$M^*=\Hom_{A}(M,A),$$ which is a $U_{-\chi}\otimes A$-module via  $x.f(m)=-f(x.m)$ for $f\in M^*$, $m\in M$, and $x\in\ggg$. It has a $X/\bbz I$-grading given by 
$$(M^*)_{\lambda+\bbz I}=\{f\in M^*\mid f(M_{\mu+\bbz I})=0 \text{ for all } \mu+\bbz I\neq-\lambda+\bbz I \}.$$
Then $M^*$ is an object in $\overline{\mathcal{C}_{\overline{A}}}$, and 
it induces a contravariant functor $(-)^*\colon\Ca\rightarrow\overline{\mathcal{C}_{\overline{A}}}$. When $M$ is free over $A$, we have $(M^*)^*\simeq M$ and $\overline{\overline{\mathcal{C}}}_{\overline{\overline{A}}}=\mathcal{C}_A$. Note that $(M^*)_{\lambda+\bbz I}\simeq (M_{-\lambda+\bbz I})^*$ for all $\lambda\in X$, we have
$$Z_A(\lambda)^*\simeq\overline{Z}_{\overline{A}}(-\lambda+2(p-1)\rho).$$
When consider the category $\mathcal{C}_A(\ggg_I)$, we also can define the dual functor $(-)^*\colon \mathcal{C}_A(\ggg_I)\rightarrow \overline{C}_{\overline{A}}(\ggg_I)$ (defined in a similar way). Denote by $\overline{Z}_{\overline{A},I}(\lambda)$ the deformed baby Verma module in $\overline{C}_{\overline{A}}(\ggg_I)$. Note that $\rho-w_i\rho$ is the sum of positive roots in $R_I$ and contained in $\bbz I$. We get similarly that
$$Z_{A,I}(\lambda)^*\simeq \overline{Z}_{\overline{A},I}(-\lambda+(p-1)(\rho-w_I\rho))\simeq\overline{Z}_{\overline{A},I}(-\lambda-\rho+w_I\rho),$$
for all $\lambda\in X$.

Let $\tau$ be the automorphism of $\ggg$, satisfies $\chi\circ\tau^{-1}=-\chi$ and induces $-w_I$ on $X$. We have $\tau(x_\alpha)=x_{-w_I\alpha}$ and $\tau(h_\alpha)=h_{-w_I\alpha}$.
Write $^\tau A$ the same $k$-algebra structure with $A$ with structural map $^\tau\pi\colon U^0\rightarrow {A}$ such that $^\tau\pi(h_\alpha)=\pi(\tau^{-1}(h))$ for $h\in U^0$.
For any $M\in\overline{\mathcal{C}}_{A}$, let $^\tau M$ denote the $U_{\chi}\otimes A$-module such that $^\tau M=M$ as an $A$-module and $x.m=\tau^{-1}(x)m$ for all $x\in U_{-\chi}$ and $m\in M$. It has a $X/\bbz I$-grading given by 
$$(^\tau M)_{\lambda+\bbz I}=M_{-\lambda+\bbz I},$$
for all $\lambda\in X$. 
Then $^\tau M$ is an object in $\mathcal{C}_A$ and $\tau$ induce a covariant functor $^\tau(-)\colon\overline{\mathcal{C}_A}\rightarrow \mathcal{C}_{^\tau A}$.
Note that $\tau(\nnn)=w^I\nnn$ and $\tau$ induces $-w_I$ on $X$, 
we have $$^\tau\overline{Z}_{A}(\lambda)\simeq Z^{w^I}_{^\tau A}(-w_I\lambda),$$
for all $\lambda\in X$. Consider the category $\overline{\mathcal{C}}_A(\ggg_I)$. Since $\tau(\ggg_I)=\ggg_I$, we can similarly define functor $^\tau(-)\colon\overline{\mathcal{C}}_A(\ggg_I)\rightarrow \mathcal{C}_{^\tau A}(\ggg_I)$. Since $\tau(\nnn_I)=\nnn_I$, we get therefore 
$$^\tau\overline{Z}_{A,I}(\lambda)\simeq Z_{^\tau A,I}(-w_I\lambda)$$
for all $\lambda\in X$.

Write $\bbd A$ the same $k$-algebra structure with $A$ but structural map $\bbd\pi\colon U^0\rightarrow  A$ extended from $h\mapsto -\pi(\tau^{-1}(h))$ for $h\in\hhh$. 
We have $\bbd A={}^\tau(\overline{A})$.
For any $M\in\Ca$, let $\bbd(M)$ denote the $U_\chi\otimes A$-module $\tau(M^*)$ with the action $x.f(m)=f(-\tau^{-1}x.m)$ for $x\in U_\chi$ and $m\in M$. 
It keeps the $X/\bbz I$-grading, which means $\bbd M_{\lambda+\bbz I}=(\bbd M)_{\lambda+\bbz I}$ as an $A$-module. We get therefore a functor $\bbd={}^\tau(-)\circ(-)^*\colon\mathcal{C}_A\rightarrow \mathcal{C}_{\bbd A}$.
If $A=F$ is a field, then $\bbd$ is an anti-equivalence functor, see \cite[Proposition 5.13]{We22}. We have
	$$\Ext^i_{\mathcal{C}_F}(M,N)\simeq \Ext^i_{\mathcal{C}_F}(\bbd N,\bbd M)$$
for all $i>0$ and $M, N\in\mathcal{C}_F$.
	
Consider the category $\Ca(\ggg_I)$. We can similarly define  $\bbd\colon\mathcal{C}_A(\ggg_I)\rightarrow \mathcal{C}_{\bbd A}(\ggg_I)$.
We have the following results.
\begin{prop} \label{duality of Z and L}
Keep the notations above.
\begin{enumerate}
\item $\bbd Z_{A,I}(\lambda)\simeq Z_{\bbd A,I}(w_I\cdot\lambda)$,
\item  $\bbd Z_A(\lambda)\simeq Z_{\bbd A}^{w^I}((w_I\cdot\lambda)^{w^I})$, and $\bbd Z^{w^I}_A(\lambda^{w^I})\simeq Z_{\bbd A}(w_I\cdot\lambda)$
\item When $A=F$ is a field. We have $\bbd L_F(\lambda)\simeq L_{\bbd F}(w_I\cdot\lambda)$.
\end{enumerate}
\end{prop}
\begin{proof}
(1)We have
$$\bbd Z_{A,I}(\lambda)={}^\tau(Z_{A,I}(\lambda)^*)={}^\tau(\overline{Z}_{\overline{A},I}(-\lambda-\rho+w_I\rho))=Z_{\bbd A,I}(w_I\cdot\lambda)$$
for all $\lambda\in X$ as desired.
 
(2)Note that $w_I\rho+w^I\rho=0$, we have
\begin{align*}
w_I(\lambda-2(p-1)\rho)+p\bbz I&=w_I\cdot\lambda-(w_I\rho-\rho)-2(p-1)w_I\rho+p\bbz I\\
&=w_I\cdot\lambda-(p-1)(\rho+w_I\rho)+p\bbz I\\
&=(w_I\cdot\lambda)^{w^I}+p\bbz I.
\end{align*} 
It follows that 
$$\bbd Z_A(\lambda)={}^\tau(Z_A(\lambda)^*)\simeq{}^\tau(\overline{Z}_{\overline{A}}(-\lambda+2(p-1)\rho))\simeq Z_{\bbd A}((w_I\cdot\lambda)^{w^I}).$$ 
Note that $w^Iw^I(\nnn)=\nnn$ and $w^I\rho$ is the sum of roots in $w^I\nnn$, we can similarly get
$$Z^{w^I}_A(\lambda)^*\simeq \overline{Z}^{w^I}_{\overline{A}}(-\lambda+2(p-1)w^I\rho)\text{ and } ^\tau\overline{Z}^{w^I}_{A}(\lambda)\simeq Z_{^\tau A}(-w_I\lambda).$$ 
Since 
$$w_I(\lambda+\rho+w^I\rho)=w_I\cdot\lambda-(w_I\rho-\rho)+w_I(\rho-w_I\rho)=w_I\cdot\lambda,$$
it follows that
$$\bbd Z_A^{w^I}(\lambda^{w^I})\simeq{}^\tau(\overline{Z}^{w^I}_{\overline{A}}(-\lambda+(p-1)(\rho+w^I\rho)))\simeq Z_{\bbd A}(w_I\cdot\lambda).$$

(3)
Note that $L_{F}(\lambda)_{\lambda+\bbz I}=Z_{F,I}(\lambda)$, we have $\bbd (L_{F}(\lambda)_{\lambda+\bbz I})=(\bbd L_F(\lambda))_{\lambda+\bbz I}$ as a $U_I\otimes F$-module. 
That is $(\bbd L_F(\lambda))_{\lambda+\bbz I}\simeq Z_{\bbd F,I}(w_I\cdot\lambda)$. Since $\bbd$ send simple object to simple object, then we have $\bbd L_{F}(\lambda)=L_{\bbd F}(w_I\cdot\lambda)$.
\end{proof}
\begin{rem}
One can compare Proposition \ref{duality of Z and L}(3) and Proposition \ref{prop long homo} to see $\bbd$ send the head of $Z_F(\lambda)$ to the scole of $\bbd Z^{w^I}_{F}(\lambda^{w^I})$ when $A=F$ is a field.
\end{rem}

\begin{lem}\label{hom ext of Z and Z^w}
	Suppose that $B$-algebra $A=F$ is a field, and $R^+_\pi\cap R_I$ equal to $\emptyset$, $\{\alpha\}$, or $R^+_I$. We have 	
\begin{enumerate}
	\item $\Hom_{\mathcal{C}_F}(Z_F(\lambda),Z_F^{w^I}(\mu^{w^I}))=
	\begin{cases}
		F, & \text{if }\lambda\in W_{I,\pi}\cdot\mu+p\bbz I,\\
		0, & otherwise.
	\end{cases}$
	\item If
	$\Ext_{\mathcal{C}_F}^1(Z_F(\lambda),Z_F^{w^I}(\mu^{w^I}))\neq0$, then $\lambda\in W_{I,\pi}\cdot\mu+p\bbz I$ with $R^+_\pi\cap R_I$ is $R^+_I$ or $\{\alpha\}$ with $\alpha\in I$.
\end{enumerate}

\end{lem}
\begin{proof}
	(1)	If $\Hom_{\mathcal{C}_F}(Z_F(\lambda),Z_F^{w^I}(\mu^{w^I}))\neq0$, then we have $Z^{w^I}_F(\mu^{w^I})_{\lambda+\bbz I}\neq0$. It follows that $\lambda+\bbz I\leq\mu+\bbz I$. Note that
	$$\Hom_{\mathcal{C}_{\bbd F}}(Z_{\bbd F}(w_I\cdot\mu),Z^{w^I}_{\bbd F}((w_I\cdot\lambda)^{w^I}))\simeq \Hom_{\mathcal{C}_F}(Z_F(\lambda),Z_F^{w^I}(\mu^{w^I})\neq0,$$
	we also can get $\mu+\bbz I\leq \lambda+\bbz I$. We get therefore $\lambda+\bbz I=\mu+\bbz I$.
    Then 
	$Z_F^{w^I}(\mu^{w^I})_{\lambda+\bbz I}=Z_F^{w^I}(\mu^{w^I})_{\mu+\bbz I}=Z_{F,I}(\mu)$. 
	The nonzero map maps $Z_F(\lambda)_{\lambda+\bbz I}=Z_{F,I}(\lambda)$ onto $Z_{F,I}(\mu)$ by simplicity.
	 Thus $\lambda\in W_{I,\pi}\cdot\mu+p\bbz I$ by Lemma \ref{propertty of g_I}.
	
	(2)Suppose $\Ext_{\mathcal{C}_F}^1(Z_A(\lambda),Z_A^{w^I}(\mu^{w^I}))\neq0$. By Lemma \ref{lem of Ext}(1), we have $\lambda+\bbz I\leq \mu+\bbz I$. 
	Note that
	$$\Ext^1_{\mathcal{C}_{\bbd F}}(Z_{\bbd F}(w_I\cdot\mu),Z_{\bbd F}(w_I\cdot\lambda)^{w^I})\simeq\Ext_{\mathcal{C}_F}^1(Z_F(\lambda),Z_F^{w^I}(\mu^{w^I}))\neq0,$$
	Similar to (1), it yields that $\lambda+\bbz I=\mu+\bbz I$. We get therefore $Z^{w^I}_F(\mu^{w^I})_{\lambda+\bbz I}=L_F(\mu)_{\lambda+\bbz I}\simeq Z_{F,I}(\mu)$. 
	If we restrict extension to the subspace of degree $\lambda+\bbz I$, it yields that
	$$\Ext^1(Z_{F,I}(\lambda),Z_{F,I}(\mu))\neq0.$$
	So the result follows from Lemma \ref{propertty of g_I}.
\end{proof}
If $M$ has a $Z$-filtration, denote by $(M\colon Z_A(\lambda))$ the multiplicity of $Z_A(\lambda)$ in the $Z$-filtration of $M$. 
Let $A=F$ be a field, recall that $Q_F(\lambda)$ is the projective cover of $L_F(\lambda)$, it has a $Z$-filtration by Lemma \ref{lem proj Z-fil}.
\begin{coro}\label{PZ ZL}
	Suppose $A=F$ is a field and $R_\pi^+\cap R_I=\emptyset$ or $R^+_\pi=\{\alpha\}$ with $\alpha\notin R_I$. Then we have
	$$(Q_F(\lambda)\colon Z_F(\mu))=[Z_F(\mu)\colon L_F(\lambda)].$$
\end{coro}
\begin{proof}
	By Lemma \ref{hom ext of Z and Z^w}, we have 	
	
	$$\Hom_{\mathcal{C}_F}(Z_F(\lambda),Z_F^{w^I}(\mu^{w^I}))=
	\begin{cases}
		F, & \lambda\in\mu+p\bbz I,\\
		0, & otherwise.
	\end{cases}$$
	and
	$$\Ext_{\mathcal{C}_F}^1(Z_F(\lambda),Z_F^{w^I}(\mu^{w^I}))=0$$
	for any $\lambda,\mu\in X$. It follows that the multiplicity of $Z_F(\lambda)$ in $Q_F(\lambda)$ is equal to the dimension of $\Hom(Q_F(\lambda),Z_F^{w^I}(\mu^{w^I}))$. Note that $[Z_F(\lambda)]=[Z_F^{w^I}(\lambda^{w^I})]$.
	Then we have 
	$$\begin{aligned}
		(Q_F(\lambda)\colon Z_F(\mu))=&\dim\Hom(Q_F(\lambda),Z_F^{w^I}(\mu^{w^I}))\\
		=&[Z_F^{w^I}(\mu^{w^I})\colon L_F(\lambda)]\\
		=&[Z_F(\mu)\colon L_F(\lambda)]
	\end{aligned}.$$
	as desired.
\end{proof}

%
%
\subsection{A-block} 

%

We are going to introduce $A$-block in the category $\Ca$.
We determine a partition of $X$ into a so-called $A$-block. We define this to be the finest possible partition such that $\lambda$ and $\mu$ belong to the same $A$-block if $\Hom_{\Ca}(Z_A(\lambda), Z_A(\mu))$ or $\Ext^1_{\Ca}(Z_A(\lambda), Z_A(\mu))$ do not vanish. This notation extends from \cite[6.9]{AJS94} in case of $I=\emptyset$.

Denote $\mathcal{D}_A$ the full subcategory of $\Ca$ containing exactly all objects with a $Z$-filtration. If $b$ is an $A$-block, denote $\mathcal{D}_A(b)$ the full subcategory of all $N$ in $\mathcal{D}_A$ where the factors in the $Z$-filtration involve only $Z_A(\lambda)$ with $\lambda\in b$, and denote $\mathcal{C}_A(b)$ the full subcategory of all modules $M$ in $\Ca$ such that $M$ is a homomorphic image of some module in $\mathcal{D}_A(b)$.

\begin{prop} Keep the notations above.
	\begin{enumerate}
		\item 	$\mathcal{D}_A=\bigoplus_{b}\mathcal{D}_A(b)$ with $\Hom(M,M')=0$ for all $M$ and $M'$ in different subcategory.
		\item   $\Ca=\bigoplus_b\Ca(b)$ with $\Hom(M,M')=0$ for all $M$ and $M'$ in different subcategory.
		\item Let $A=F$ be a field, then $F$-block concides with usual one$\colon$ the finest partition of $X$ such that $\lambda$ and $\mu$ lie in the same block if $L_F(\lambda)$ and $L_F(\mu)$ have a nontrival extension.
	\end{enumerate}	
\end{prop}
\begin{proof}
	(1) It follows directly from the definition of $A$-block.
	
	(2) Suppose $M$ and $M'$ are in $\mathcal{C}_A(b)$ and  $f\in\Hom_{\Ca}(M,N)$. By Lemma \ref{lem proj Z-fil}, there exists epimorphisms $p_M\colon Q\twoheadrightarrow M$ and $p_{M'}\colon Q'\twoheadrightarrow M'$ with $Q$ and $Q'$ are projective objects with $Z$-filtration. We may assume that $Q$ and $Q'$ are both in $\mathcal{D}_A(b)$ by (1). Since $Q$ is projective, there exists a homomorphism $f'\colon Q\rightarrow Q'$ such that $f\circ p_M=p_{M'}\circ f'$, then  $f=0$ whenever $f'=0$. 
	Note that $M=p_M(Q)=\bigoplus_b p_M(Q_b)$, where $Q_b\in\mathcal{D}_A(b)$. Then we have $p_M(Q_b)\in\mathcal{C}_A(b)$ and hence $\mathcal{C}_A=\bigoplus_b\Ca(b)$.
	
	(3) If $\lambda$ and $\lambda'$ belong to $F$-block $b$ and $b'$, we have $L_F(\lambda)\in\mathcal{C}_F(b)$ and $L_F(\lambda')\in\mathcal{C}_F(b')$. When $b\neq b'$, then $\Ext^1(L_F(\lambda),L_F(\mu))=0$. It shows that each block is contained in a $F$-block.
	
	For every block $B$, denote $\mathcal{C}_F(B)$ the full subcategory of $\mathcal{C}_F$ with objects has composition factors only involves $L_F(\lambda)$ with $\lambda\in B$. Suppose now $\lambda\in B$. Since $Z_F(\lambda)$ is indecomposable, it belong to $\mathcal{C}_F(B)$. If $\lambda'$ belong to a different block $B'$, then $\Hom_{\mathcal{C}_F}(Z_F(\lambda),Z_F(\lambda'))=0=\Ext_{\mathcal{C}_F}(Z_F(\lambda),Z_F(\lambda'))$ by the property of block. Hence, each $F$-block is contained in a block. Thus block and $F$-block coincide.
\end{proof}

We want to determine the $A$-block when $A$ is a field at first.
Recall that $B=S[(h_\alpha^{p-1}-1)^{-1}\mid\alpha\in R]$ and definition of $n_\beta(\lambda)$ in \ref{sec twist}. 
When $A$ is a $B$-algebra, then we have
  $$\langle\lambda+\rho,\beta^\vee\rangle\equiv n_\beta(\lambda)\;(\mod p),$$
 for each $\lambda\in X$ and $\beta\in R_\pi$. 

Let $W_{\pi,p}=\langle s_{\alpha,mp}\mid \alpha\in (R_\pi\cap I)\cup (R_\pi\backslash R_I), m\in\bbz\rangle$. For any $M\in\mathcal{C}_F$, denote $[M\colon L_F(\lambda)]$ the multiplicity of $L_F(\lambda)$ as a composition factor of $M$. 
we first investigate the case of field.
 
\begin{lem} \label{lem block field}
	Suppose $B$-algebra $A=F$ is a field.
	\begin{enumerate} 
		\item 	Let $\lambda, \mu\in X$ with $\lambda+\bbz I\neq\mu+\bbz I$. Suppose $[Z_F(\lambda)\colon L_F(\mu)]\neq0$, then  $[Z_F(\lambda-n_\beta(\lambda)\beta)\colon L_F(\mu)]\neq0$ for some $\beta\in R_{\pi}\backslash R_I$ with $n_\beta(\lambda)\neq 0$.
		\item Let $\lambda, \mu\in X$ with $\lambda+\bbz I=\mu+\bbz I$ and suppose $R_\pi\cap R_I=\emptyset,\{\pm\alpha\}$ or $R_I$. If  $[Z_F(\lambda)\colon L_F(\mu)]\neq0$ , then
		$\mu\in W_{I,\pi}\cdot\lambda+p\bbz I$.
		\item  Each $F$-block of $\lambda\in X$ belong to  $W_{\pi,p}\cdot\lambda+p\bbz I$.
	\end{enumerate}
\end{lem}
\begin{proof}
	(1)Since $\lambda+\bbz I\neq\mu+\bbz I$. Then $[Z_F(\lambda)\colon L_F(\mu)]\neq0$ yields that there exists $w\in W^I$ and simple root $\alpha$ with $w\alpha\in R_\pi\backslash R_I$ such that $L_F(\mu)$ lies in the kernel of 
	$\varphi\colon Z^w_F(\lambda^{w})\rightarrow Z^{s_\alpha w}_F (\lambda^{w}-(p-1)w\alpha)$. Let $\beta=w\alpha$,
	then $L_F(\mu)$ is in the image of $Z_F^{w}(\lambda^{w}-n_\beta(\lambda)\beta)$ by Lemma \ref{twist homo property}(2) and
	$$\langle\lambda+\rho,\beta^\vee\rangle=\langle\lambda^{w},\beta^\vee\rangle+\langle\rho+(p-1)(\rho-w\rho),w\alpha^\vee\rangle\equiv\langle\lambda^w,\beta^\vee\rangle+1 \;(\mod p).$$
	It follows that
	$$[Z_F(\lambda-n_\beta(\lambda)\beta):L_F(\mu)]=[Z^w_F(\lambda^w-n_\beta(\lambda)\beta):L_F(\mu)]\neq0.$$

	(2)Note that $Z_F(\lambda)_{\lambda+\bbz I}=Z_{F,I}(\lambda)$ and $L_F(\mu)_{\lambda+\bbz I}=Z_{F,I}(\mu)$ are irreducible. $[Z_F(\lambda)\colon L_F(\mu)]\neq0$ if and only if $Z_{F,I}(\lambda)\simeq Z_{F,I}(\mu).$ Then the results follows from Lemma \ref{propertty of g_I}.
	
	(3)Suppose that $\Ext^1_{\mathcal{C}_F}(L_F(\lambda),L_F(\mu))\neq 0$ with $\lambda,\mu\in X$.
	
	If $\lambda+\bbz I\nleq\mu+\bbz I$, it follow from
	 Lemma \ref{lem of Ext} that $$\Ext^1_{\mathcal{C}_F}(L_F(\lambda),L_F(\mu))=\Hom_{\mathcal{C}_F}(\Rad Z_F(\lambda),L_F(\mu))\neq0.$$
	Note that $\lambda-n_\beta(\lambda)\beta=s_{\beta,mp}\cdot\lambda$ where $mp=\langle\lambda+\rho,\beta^\vee\rangle-n_\beta(\lambda)$.
	By repeatedly applying (1), there exists $\lambda'\in W_{\pi,p}\cdot\lambda+p\bbz I$ such that $[Z_F(\lambda')\colon L_F(\mu)]\neq0$ and $\lambda'+\bbz I=\mu+\bbz I$. Then the result follows by (2) since $W_{I,\pi}\subset W_{\pi,p}$.
	
	If $\lambda+\bbz I=\mu+\bbz I$, we have 
	$L_F(\lambda)_{\lambda+\bbz I}=Z_{F,I}(\lambda)$ and $L_F(\mu)_{\lambda+\bbz I}=Z_{F,I}(\mu)$. Similar to
	Proposition \ref{lem of Ext}(2), we could get
	$\Ext^1(Z_{F,I}(\lambda),Z_{F,I}(\mu))\neq0$.
	Then $\mu\in W_{I,\pi}\cdot\lambda+p\bbz I\subset W_{\pi,p}\cdot\lambda+p\bbz I$ as desired. 
	
	If $\lambda+\bbz I<\mu+\bbz I$, we have	
	\begin{align*}
		\Ext^1_{\mathcal{C}_F}(L_F(\lambda),L_F(\mu))&\simeq \Ext^1_{\mathcal{C}_{\mathbb{D}F}}(\bbd L_{F}(\mu),\bbd L_{ F}(\lambda))\\
		&=\Ext^1_{\mathcal{C}_{\mathbb{D}F}}(L_{\bbd F}(w_I\cdot\mu),L_{\bbd F}(w_I\cdot\lambda))\\
		&\simeq\Hom_{\bbd F}(\Rad Z_{\bbd F}(w_I\cdot\mu),L_{\bbd F}(w_I\cdot\lambda)).
	\end{align*}
	Note that $R_{\bbd\pi}=w_I R_\pi$ and $\langle w_I\cdot\mu+\rho,w_I\alpha^\vee\rangle=\langle\mu+\rho,\alpha^\vee\rangle$. 
	Similar to the case of $\lambda+\bbz I\nleq\mu+\bbz I$, we have
	 $w_I\cdot\lambda\in W_{\bbd\pi,p}\cdot(w_I\cdot\mu)+p\bbz I$. 
	Note that $w_I(R^+_I)=-R^+_I$ and $w_I(I)=-I$. We have 
	$$W_{\bbd\pi,p}=\langle s_{\alpha,mp}\mid \alpha\in w_I((R_\pi\cap I)\cup (R_\pi\backslash R_I)), m\in\bbz\rangle.$$	
	By conjugation in the Weyl group, we conclude that $\mu\in W_{\pi,p}\cdot\lambda+p\bbz I$ as desired.

\end{proof}

\begin{prop}\label{A-block}
Let $\lambda\in X$ and suppose $A$ is a $B$-algebra, then each $A$-block of $\lambda$ belong to $W_{\pi,p}\cdot\lambda+p\bbz I$.
\end{prop}
\begin{proof}
  For any $A$-algebra $A'$, let $$E_{A'}^i=\Ext_{\mathcal{C}_{A'}}^i(Z_{A'}(\lambda),Z_{A'}(\mu)).$$
Suppose now $\lambda,\mu\in X$ such that $E^1_A\neq0$ or $E^0_A\neq0$.
	
	If $E_{A/\mmm}^1\neq0$ for some maximal ideal $\mmm$ of $A$, we have $\lambda\in W_{\pi,p}\cdot\mu$ since $R_\pi$ for $A/\mmm$ is contained in $R_\pi$ for $A$.
	
	If $E_{A/\mmm}^1=0$ for all maximal ideal $\mmm$ of $A$, we have $E_{A/\mmm}^i=0$ for all maximal ideal $\mmm$ of $A$ with $i>0$. By Lemma \ref{lem ext base change}, we have $E^i_A=0$ for all $i>0$, then $E_A^0\neq0$. Since
	$E^0_A\otimes A/\mmm\simeq E^0_{A/\mmm}$
	  and Nakayama lemma,
     there exists a maximal ideal $\mmm$ such that $E_{A/\mmm}^0\neq0$, which yields that $\lambda\in W_{\pi,p}\cdot\mu+p\bbz I$ by Lemma \ref{lem block field}.
\end{proof}
\begin{coro}\label{irr and proj of Z}
	Suppose $B$-algebra $A=F$ is a field and let $\lambda\in X$. If $n_\beta(\lambda)=0$ for all $\beta\in R_\pi\backslash R_I$ and $ R_\pi\cap I=\emptyset$, then all $Z_F(\lambda)$ are projective and irreducible.
\end{coro}
\begin{proof}
	The irreducibility of $Z_F(\lambda)$ comes from Corollary \ref{coro irreducible of Z}. By Proposition \ref{A-block}, the $F$-block of each $\lambda\in X$ is equal to $\lambda+p\bbz I$. When $\mu\notin\lambda+p\bbz I$, we have
	$$\Ext^1_{\mathcal{C}_F}(L_F(\lambda),L_F(\mu))=\Ext^1_{\mathcal{C}_F}(Z_F(\lambda),Z_F(\mu))=0.$$
	When $\mu\in\lambda+p\bbz I$. By Lemma \ref{propertty of g_I}, we have $\Ext^1_{\mathcal{C}_F(\ggg_I)}(Z_{F,I}(\lambda),Z_{F,I}(\mu))=0$. We deduce from Lemma \ref{lem of Ext} that $$\Ext^1_{\mathcal{C}_F}(L_F(\lambda),L_F(\mu))=\Ext^1_{\mathcal{C}_F}(Z_F(\lambda),Z_F(\mu))=0.$$
	
	We conclude that $\Ext^1(L_F(\lambda),M)=0$ for all simple module $M$ and $\lambda\in X$. Therefore $Z_F(\lambda)=L_F(\lambda)$ is projective for all $\lambda\in X$.
\end{proof}
\begin{coro}\label{Z are progenerator}
	Suppose $A$ is a $B^\emptyset$-algebra or $B^\alpha$-algebra with $\alpha\in R^+_I\backslash I$. For any $\lambda\in X$, all $Z_A(\lambda)$ are projective and generate $\Ca$.
\end{coro}
\begin{proof}
	For any any maximal ideal $\mmm$ of $A$, 
    the analogue of $R_\pi$ for $A/\mmm$ contained in $R_\pi$ for $A$.
	By Lemma \ref{lem R pi} and Corollary \ref{irr and proj of Z}, all $Z_{A/\mmm}(\lambda)$ are projective in $\mathcal{C}_{A/\mmm}$. We can now deduce from Corollary \ref{coro projective from field to ring} that $Z_A(\lambda)$ is projective in $\Ca$ for each $\lambda\in X$.
	Note that the $A$-block of $\lambda$ is $\lambda+p\bbz I$.
	 Then any $M\in\mathcal{C}(\lambda+p\bbz I)$ is a homomorphic image of $Z_A(\lambda)^{\oplus m}$ for some $m>0$. So all $Z_A(\lambda)$ generate $\Ca$.
\end{proof}

\subsection{Extensions in rank 1 cases}\label{subsec proj}

We focus on the projective object in $\Ca$. We shall always assume $A$ and any $A$-algebra are Noetherian and commutative.
We are not going to determine a general projective object but under some restrictions on $A$ or $\pi$.

In this section, we keep the assumption that $A$ is an integral domain and $R^+_\pi=\{\alpha\}$ with $\alpha\notin R_I^+\backslash I$. Let $\lambda\in X$ and $n\in\bbz$ with $0\leq n<p$ such that  $\langle\lambda+\rho,\alpha^\vee\rangle\equiv p-n \;(\mod p)$. Put $\alpha\uparrow\lambda=\lambda+n\alpha$.

\begin{prop}
Suppose $A$ is a $B^\alpha$-algebra with $\alpha\notin R^+_I\backslash I$. If $\alpha\uparrow\lambda=\lambda$, then $Z_A(\lambda)$ is projective.
\end{prop}
\begin{proof}
By Corollary \ref{coro projective from field to ring}, we only need to show $Z_{A/\mmm}(\lambda)$ is projective for all maximal ideal $\mmm$ of $A$. If $\pi(h_\alpha)\notin\mmm$, then the analogue of $R_\pi$ for $A/\mmm$ is equal to $\emptyset$. We have $Z_{A/\mmm}(\lambda)$ is projective by Corollary \ref{irr and proj of Z}. If $\pi(h_\alpha)\in\mmm$, then the analogue of $R_\pi$ for $A/\mmm$ is equal to $\{\pm\alpha\}$. We deduce from $\alpha\uparrow\lambda=\lambda$ that $|W_{I,\pi}\cdot\lambda|=1$. Then $Z_{A/\mmm,I}(\lambda)$ is projective by Lemma \ref{propertty of g_I}, and then
$$\Ext^1_{\mathcal{C}_{A/\mmm}(\ggg_I)}(Z_{A/\mmm,I}(\lambda),Z_{A/\mmm,I}(\mu))=0$$
for each $\mu\in X$.
Similar to Lemma \ref{lem of Ext}(2), we could get
$$\Ext^1_{\mathcal{C}_{A/\mmm}}(Z_{A/\mmm}(\lambda),L_{A/\mmm}(\mu))=0$$ for any irreducible module $L_{A/\mmm}(\mu)$.  Hence, $Z_{A/\mmm}(\lambda)$ is projective.
\end{proof}

\begin{lem}\label{projective in rank 1 case}
	Suppose $B^\alpha$-algebra $A=F$ is a field with $\alpha\notin R^+_I\backslash I$.
 If $\alpha\uparrow\lambda\neq\lambda$, then $Q_F(\lambda)$ has a $Z$-filtration such that the all composition factors are $Z_F(\lambda)$ and $Z_F(\alpha\uparrow\lambda)$ and each occurs once. In particular, we have
	$$\Ext^1_{\mathcal{C}_F}(Z_{F}(\lambda),Z_{F}(\alpha\uparrow\lambda))=F.$$
\end{lem}
\begin{proof}
In case $\alpha\in I$. We have $Z_F(\lambda)$ is irreducible for each $\lambda\in X$ by Corollary \ref{coro irreducible of Z} and $Z_F(\lambda)\simeq Z_F(\mu)$ if and only if $\mu\in\lambda+p\bbz I$ or $\mu\in\alpha\uparrow\lambda+p\bbz I$ by Corollary \ref{coro iso class of Z}.
By \cite[Lemma 10.9]{Ja98}, we have 
$$\dim Q_{F}(\lambda)=p^{\dim\nnn}\displaystyle\sum_{\mu+pX\in X/pX}[Z_{F}(d\mu)\colon L_{F}(d\lambda)]=2p^{\dim\nnn}.$$
Note that the $F$-block of $\lambda$ is $(\lambda+p\bbz I)\cup(\alpha\uparrow\lambda+p\bbz I)$.
It follows that $Q_F(\lambda)$ has a $Z$-filtration such that the only two factors are both isomorphic to $Z_F(\lambda)$.
Then we have a projective resolutin of $Z_{F}(\lambda)$
$\colon$
$$\cdots\xlongrightarrow{\bar{p}} Q_{F}(\lambda)\xlongrightarrow{\bar{p}}Q_{F}(\lambda)\xlongrightarrow{p} Z_{F}(\lambda)\longrightarrow 0,$$
where $p\colon Q_F(\lambda)\rightarrow Z_F(\lambda)$ is the projection and we also denote $Q_F(\lambda)\xrightarrow{p} Z_F(\lambda)\hookrightarrow Q_F(\lambda)$ by $\bar{p}$. By applying the functor $\Hom_{\mathcal{C}_F}(-,Z_F(\lambda))$, we get the complex
$$0\rightarrow\Hom_{\mathcal{C}_F}(Z_F(\lambda),Z_F(\lambda))\rightarrow\Hom_{\mathcal{C}_F}(Q_F(\lambda),Z_F(\lambda))\rightarrow\Hom_{\mathcal{C}_F}(Q_F(\lambda),Z_F(\lambda))\rightarrow\cdots$$
 Since
  ${p}\bar{p}=0$, we get theroefore
$$\Ext^1_{\mathcal{C}_F}(Z_F(\lambda),Z_F(\lambda))\simeq\Hom_{\mathcal{C}_F}(Q_F(\lambda),Z_F(\lambda))=F.$$

In the case of $\alpha\in R^+\backslash R^+_I$.
We now use the notations in Proposition \ref{prop long homo}. Suppose $\alpha=w_i\alpha_i$, then we have $\varphi_j$ are all isomorphism with $j\neq i$ and there exists $\psi\colon Z^{w_i}_F(\lambda^{w^I}-(p-n)\alpha)\rightarrow Z_F^{w_i}(\lambda)$ such that $\ker\varphi_i\simeq\im\psi$ by Proposition \ref{twist homo property}(2).
Note that $$\dim L_F(\lambda)=\dim\im\varphi=\dim\im\varphi_i=p^{\dim\nnn-1}(p-n),$$ which yields that $\dim L_F(\lambda-(p-n)\alpha)=p^{\dim\nnn-1}n$. 
Also, we have 
$$\dim\im\psi=\dim\ker\varphi_i'=p^{\dim\nnn-1}n.$$ We get therefore $L_F(\lambda-(p-n)\alpha)\simeq\im\psi\simeq\ker\varphi_i\simeq\ker\varphi$. So we have an exact sequence
$$0\rightarrow L_F(\lambda-(p-n)\alpha)\rightarrow Z_F(\lambda)\rightarrow L_F(\lambda)\rightarrow0.$$  
It yields by Corollary \ref{PZ ZL} that an exact sequence
$$0\rightarrow Z_F(\alpha\uparrow\lambda)\longrightarrow Q_F(\lambda)\longrightarrow Z_F(\lambda)\rightarrow0.$$
Then we get a projective resolution of $Z_F(\lambda)$
$$\cdots\xrightarrow{f} Q_F(\lambda+p\alpha)\xrightarrow{g} Q_F(\alpha\uparrow\lambda)\rightarrow Q_F(\lambda)\rightarrow Z_F(\lambda)\rightarrow0.$$ 
By applying the functor $\Hom_{\mathcal{C}_F}(-,Z_F(\alpha\uparrow\lambda))$, we get the complex
$$0\rightarrow\Hom_{\mathcal{C}_F}(Z_F(\lambda),Z_F(\alpha\uparrow\lambda))\rightarrow\Hom_{\mathcal{C}_F}(Q_F(\lambda),Z_F(\alpha\uparrow\lambda))\rightarrow\Hom_{\mathcal{C}_F}(Q_F(\alpha\uparrow\lambda),Z_F(\alpha\uparrow\lambda))\rightarrow0.$$
So we have $$\Ext^1(Z_F(\lambda),Z_F(\alpha\uparrow\lambda))\simeq\Hom_{\mathcal{C}_F}(Q_F(\alpha\uparrow\lambda),Z_F(\alpha\uparrow\lambda))=F,$$
as desired.
\end{proof}
 Let $K$ be the fraction field of integral domain $A$. Suppose $R_\pi^+=\{\alpha\}$ with $\alpha$ being simple.
Denote the generator $v_0\in Z_A(\lambda)$ and $v_1\in Z_A(\lambda+n\alpha)$. For any $b\in K$, let $y_b=v_0+E_{-\alpha}^{(n)}v_1\in Z_K(\lambda)\oplus Z_K(\lambda+n\alpha)$ and
$$Y(b)=U^-v_0A+U^-y_bA\subset Z_K(\lambda)\oplus Z_K(\lambda+n\alpha).$$

We study in the case of $\alpha$ being simple at first.
We have the following results from \cite[8.1-8.4]{AJS94}.
\begin{lem}\label{lem of ext of rank 1}
	Suppose $A$ is a $B^\alpha$-algebra with simple root $\alpha$ and $\alpha\uparrow\lambda\neq\lambda$. Let $b\in K$. 
	\begin{enumerate}
		\item $Y(b)$ is an object in $\Ca$ if and only if $\pi(h_\alpha)b\in A$.
		\item Suppose $\pi(h_\alpha)b\in A$, then there exists an exact sequence in $\Ca$
		$$[b]\colon0\longrightarrow Z_A(\alpha\uparrow\lambda)\xlongrightarrow{f} Y(b)\xlongrightarrow{g} Z_A(\lambda)\longrightarrow0,$$
		where $f(v_1)=v_1$ and $g(y_b)=v_0$. In addition, the exact sequence split if and only if $b\in A$.
		\item We have an isomorphism
		$$\Ext^1_{\Ca}(Z_A(\lambda),Z_A(\alpha\uparrow\lambda))\simeq A\pi(h_\alpha)^{-1}/A.$$
	\end{enumerate}
\end{lem}
\begin{coro}\label{projective as generators}
	Suppose $A$ is a $B^\alpha$-algebra with simple root $\alpha$ and $\alpha\uparrow\lambda\neq\lambda$. Then
	$Y(\pi(h_\alpha)^{-1})$ is projective in $\Ca$ and the exact sequence $[\pi(h_\alpha)^{-1}]$ generates $\Ext^1_{\Ca}(Z_A(\lambda),Z_A(\alpha\uparrow\lambda))$.
\end{coro}
\begin{proof}
	Let $\mmm$ be a maximal ideal of $A$. If $\pi(h_\alpha)\notin\mmm$, then the analogue of $R_\pi$ for $A/\mmm$ is emptyset. We have $Z_{A/\mmm}(\lambda)$ and $Z_{A/\mmm}(\alpha\uparrow\lambda)$ are projective by Corollary \ref{irr and proj of Z}, so is $Y(\pi(h_\alpha)^{-1})\otimes A/\mmm$.
	If $\pi(h_\alpha)\in\mmm$, it follows from  Lemma \ref{projective in rank 1 case} and Lemma \ref{lem of ext of rank 1} that $Y(\pi(h_\alpha)^{-1})\otimes A/\mmm$ is non-split and isomorphic to $Q_{A/\mmm}(\lambda)$.
	We conclude that $Y(\pi(h_\alpha)^{-1})\otimes A/\mmm$ is projective for each maximal ideal $\mmm$, so is $Y(\pi(h_\alpha)^{-1})$.  The remaining part comes from Lemma \ref{lem of ext of rank 1}(3).
\end{proof}

Now we back to the general case.
\begin{prop}\label{prop ext of rank 1}
Suppose $A$ is a $B^\alpha$-algebra and $\alpha\uparrow\lambda\neq\lambda$. 
\begin{enumerate}
	\item $\Ext^1_{\Ca}(Z_A(\lambda),Z_A(\alpha\uparrow\lambda))\simeq A\pi(h_\alpha)^{-1}/A.$	
	\item There exists a projective objective $Q$ in $\mathcal{C}_A$ with an exact sequence
	$$0\rightarrow Z_A(\alpha\uparrow\lambda)\longrightarrow Q\longrightarrow Z_A(\lambda)\rightarrow0.$$ 
\end{enumerate}
\end{prop}
\begin{proof}
We only need proof the case of $\alpha\in R^+\backslash R^+_I$ not simple. Note that there exists $w\in W^I$ such that $\alpha$ is simple  with respect to roots in $w\nnn$. Note that 
$\langle\lambda^w+w\rho,\alpha^\vee\rangle=\langle\lambda+\rho,\alpha^\vee\rangle\equiv p-n\; (\mod p)$, we have $\alpha\uparrow\lambda^w=\lambda^w+n\alpha$.
Similar to Lemma \ref{lem of ext of rank 1} and Corollary \ref{projective as generators}, there exists projective object $Q$ with an exact sequence
$$0\rightarrow Z^w_A(\alpha\uparrow\lambda^w)\longrightarrow Q\longrightarrow Z^w_A(\lambda^w)\rightarrow0.$$
Then (2) followed by isomorphism $Z^w_A(\lambda^w)\simeq Z_A(\lambda)$ and $Z^w_A(\alpha\uparrow\lambda^w)\simeq Z_A(\alpha\uparrow\lambda)$.
We have $$\Ext^1_\Ca(Z_A(\lambda),Z_A(\alpha\uparrow\lambda))=\Ext^1_\Ca(Z^w_A(\lambda^w),Z^w_A(\alpha\uparrow\lambda^w))=A\pi(h_\alpha)^{-1}/A$$ as desired.
\end{proof}
%

\section{A combinatorial category}
In this section, we keep the assumption that $A=\hat{S}$ to be a complection of $S=S(\hhh)$ at the augmentation ideal of $S$. Then $\hat{S}$ is a Noetherian local integral domain over $B$ and flat over $S$. We denote $A^\emptyset=\hat{S}^\emptyset=\hat{S}[h_\beta^{-1}\mid \beta\in R^+]$ and $A^\alpha=\hat{S}^\alpha=\hat{S}[h_\beta^{-1}\mid \beta\in R^+\backslash\{\alpha\}]$ for any $\alpha\in R^+$. 
We have naturally embedding $S\hookrightarrow \hat{S}\hookrightarrow \hat{S}^\alpha\hookrightarrow \hat{S}^\emptyset$ and
the structural maps are the natural embeddings of $U^0\simeq S$. For any $\hat{S}$-module $M$, we denote $M^\emptyset=M\otimes_{\hat{S}}\hat{S}^\emptyset$ and  $M^\alpha=M\otimes_{\hat{S}}\hat{S}^\alpha$. For all $\lambda\in X$, we use abbreviations that $Z^\emptyset(\lambda)=Z_{\hat{S}^\emptyset}(\lambda)$ in $\mathcal{C}^\emptyset=\mathcal{C}_{\hat{S}^\emptyset}$ and $Z^\alpha(\lambda)=Z_{\hat{S}^\alpha}(\lambda)$ in $\mathcal{C}^\alpha=\mathcal{C}_{\hat{S}^\alpha}$

We have nice properties of the assumptions from \cite[9.1-9.2]{AJS94}.

\begin{lem}
Suppose that $p\neq2$ if $R$ has two root lengths and $p\neq3$ if $R$ has a component of type $G_2$. (The characteristic $p$ is prime to all entries of the Cartan matrix.) Then 
\begin{enumerate}
\item $\hat{S}=\bigcap_{\beta\in R^+}\hat{S}^\beta$,
\item If $M$ is a flat $\hat{S}$-module, $M=\bigcap_{\beta\in R^+}M^\beta$.
\end{enumerate}
\end{lem}


Let $\Omega$ be an $(\cdot)$-orbit of $W_p$ in $X$  under the dot action. Since $R_\pi=R$ for $\mathcal{C}_{\hat{S}}$, then $\Omega$ is a $\hat{S}$-block of $\mathcal{C}_{\hat{S}}$ by Proposition \ref{A-block}.

Fix a ($\cdot$)-orbits $\Omega$ of $W_p$ in $X$. We define a combinatorial category $\mathcal{K}=\mathcal{K}(\Omega,\hat{S})$ 
consisting of
	\begin{itemize}
		\item Objects$\colon$ $\mathcal{M}=((\mathcal{M}(\lambda))_{\lambda\in\Omega},\mathcal{M}(\lambda,\alpha)_{\lambda\in\Omega,\alpha\in R^+})$, where each $\mathcal{M}(\lambda)$ is finitely generated $\hat{S}^\emptyset$-module and each $\mathcal{M}(\lambda,\alpha)$ is finitely generated $\hat{S}^\alpha$-submodule of $\mathcal{M}(\lambda)\oplus\mathcal{M}(\alpha\uparrow\lambda)$.
		\item Morphisms$\colon$ $\phi\colon\mathcal{M}\rightarrow\mathcal{N}$ is a collection of $\hat{S}^{\emptyset}$-homomorphisms $\phi_\lambda\colon\mathcal{M}(\lambda)\rightarrow\mathcal{N}(\lambda)$ which satisfies $(\phi_{\lambda}\oplus\phi_{\alpha\uparrow\lambda})\mathcal{M}(\lambda,\alpha)\subset\mathcal{N}(\lambda,\alpha)$.
	\end{itemize}
  
In $\mathcal{C}_{\hat{S}^\alpha}$, we put $Q^\alpha(\lambda)=Z^\alpha(\lambda)$ if $\alpha\in R^+_I\backslash I$ or $\alpha\uparrow\lambda=\lambda$ and $Q^\alpha(\lambda)=Q$ being the projective object
in Corollary \ref{prop ext of rank 1} otherwise. We have all $Q^\alpha(\lambda)$ are projective and generate $\mathcal{C}^\alpha$.
 
Denote $\mathcal{FC}_{\hat{S}}(\Omega)$ be the full subcategory of $\mathcal{C}_{\hat{S}}(\Omega)$ consisting of objects in $\mathcal{C}_{\hat{S}}(\Omega)$ that are flat over $\hat{S}$. All projective object of $\mathcal{C}_{\hat{S}}(\Omega)$ are contained in $\mathcal{FC}_{\hat{S}}(\Omega)$ by Lemma \ref{lem enough proj}(2).
Then we can define a functor
$$\mathbb{V}=\mathbb{V}_{\Omega}:\mathcal{FC}_{\hat{S}}(\Omega)\longrightarrow \mathcal{K}(\Omega,\hat{S}) $$ as follows.
For $M\in\mathcal{FC}_{\hat{S}}(\Omega)$, let 
$\mathbb{V}M(\lambda)=\Hom_{\mathcal{C}^\emptyset}(Z^\emptyset(\lambda),M^\emptyset)$, 
and
$$\mathbb{V}M(\lambda,\alpha)=\im	(\Hom_{\mathcal{C}^\alpha}(Q^\alpha(\lambda),M^\alpha)\hookrightarrow\Hom_{\mathcal{C}^\alpha}(Q^\alpha(\lambda),M^\alpha)\otimes\hat{S}^\emptyset
)).$$
We have $\mathbb{V}M(\lambda)$ is finitely generated over $\hat{S}^\emptyset$ since $\hat{S}^\emptyset$ is Noetherian.
 Since $\hat{S}^\alpha$ and $\hat{S}^\emptyset$ are flat over $\hat{S}$, we deduce from Lemma \ref{lem base change}(3) that
 $$\Hom_{\mathcal{C}^\alpha}(Q^\alpha(\lambda),M^\alpha)\otimes\hat{S}^\emptyset=\Hom_{\mathcal{C}^\emptyset}(Q^\alpha(\lambda)\otimes \hat{S}^\emptyset,M^\emptyset).$$
 Note that
  $$Q^{\alpha}(\lambda)\otimes_{\hat{S}^\alpha}\hat{S}^\emptyset=
 \begin{cases}
 	Z^\emptyset(\lambda),&\text{if } \alpha\in R^+_I\backslash I \text{ or } \alpha\uparrow\lambda=\lambda,\\
 	Z^\emptyset(\lambda)\oplus Z^\emptyset(\alpha\uparrow\lambda),& \text{otherwise}.
 \end{cases}$$
It follows that $\mathbb{V}M(\lambda,\alpha)\subset\mathbb{V}M(\lambda)\oplus\mathbb{V}M(\alpha\uparrow\lambda)$, and then the  functor is well-defined.

\begin{prop}\label{prop fully faithful}
Keep the assumption on $p$. Then $\mathbb{V}$ is a fully faithful functor, that is
$$\Hom_{\mathcal{C}_{\hat{S}}(\Omega)}(M,N)\simeq\Hom_{\mathcal{K}}(\mathbb{V}M,\mathbb{V}N)$$
for any $M, N\in\mathcal{FC}_{\hat{S}}(\Omega)$
\end{prop}
\begin{proof}
Since $Z^\emptyset(\lambda)$ is a projective generator of $\mathcal{C}^\emptyset(\lambda+p\bbz I)$ (see Corollary \ref{Z are progenerator}) and $\End(Z^\emptyset(\lambda))=\hat{S}^\emptyset$, 
then the functor $$\Hom_\mathcal{C^\emptyset}(Z^\emptyset(\lambda),-)\colon  \mathcal{C}^{\emptyset}(\lambda+p\bbz I)\rightarrow \hat{S}^\emptyset\text{-mod}$$
 is fully faithful. We get therefore
$$\Hom_{\mathcal{C}^{\emptyset}}(M^\emptyset,N^\emptyset)\simeq\bigoplus_{\lambda\in\Omega}\Hom_{\hat{S}^\emptyset}(\mathbb{V}M(\lambda),\mathbb{V}N(\lambda)).$$
for any $M, N\in\mathcal{FC}_{\hat{S}}(\Omega)$. 
Note that we have a commutative diagram
$$
\xymatrix{\Hom_{\mathcal{C}_{\hat{S}}}(M,N)\ar[r]\ar[d] & \Hom_{\mathcal{K}(\Omega,\hat{S})}(\mathbb{V}M,\mathbb{V}N)\ar[d]\\
\Hom_{\mathcal{C}^{\emptyset}}(M^\emptyset,N^\emptyset)\ar[r]^(0.36){\sim} &\bigoplus_{\lambda\in\Omega}\Hom_{\hat{S}^\emptyset}(\mathbb{V}M(\lambda),\mathbb{V}N(\lambda))
}
$$
Since two vertical maps are inclusions, we have
 $\Hom_{\mathcal{C}}(M,N)\rightarrow\Hom_{\mathcal{K}}(\mathbb{V}M,\mathbb{V}N)$ is injective.

Let $\varphi=(\varphi_\lambda)_{\lambda\in\Omega}\in\Hom_{\mathcal{K}}(\mathbb{V}M,\mathbb{V}N)$, then there exists $f\in\Hom_{\mathcal{C}^{\emptyset}}(M^\emptyset,N^\emptyset)$ such that $\varphi$ is induced by $f$.  For all $\lambda\in\Omega$, the map $\varphi$ induces a map from $\Hom(Q^\alpha(\lambda),M^\alpha)$ to $\Hom(Q^\alpha(\lambda),N^\alpha)$ which sends $g$ to $f\circ g$. It follows that $f(M^\alpha)\subset N^\alpha$. Now we have $$f(M)=f(\bigcap_{\alpha}M^\alpha)=\bigcap_{\alpha}f(M^\alpha)\subset\bigcap_{\alpha}N^\alpha= N,$$
which shows the surjectivity.
\end{proof}
\begin{lem} Let $M\in\mathcal{FC}_{\hat{S}}(\Omega)$, and $N$ be the submodule of $M$ such that $N$ and $M/N$ are flat over $\hat{S}$.  
We have $\mathbb{V}(M)/\mathbb{V}(N)\simeq\mathbb{V}(M/N)$
\end{lem}
\begin{proof}
The flatness yields that $N$ and $M/N$ are also objects in $\mathcal{FC}_{\hat{S}}(\Omega)$. Since $Z^\emptyset(\lambda)$ and $Q^\alpha(\lambda)$ is projective. We have
$$0\longrightarrow\Hom_{\mathcal{C}^\emptyset}(Z^\emptyset(\lambda),N^\emptyset)\longrightarrow\Hom_{\mathcal{C}^\emptyset}(Z^\emptyset(\lambda),M^\emptyset)\longrightarrow\Hom_{\mathcal{C}^\emptyset}(Z^\emptyset(\lambda),(M/N)^\emptyset)\longrightarrow0,$$
and
$$0\longrightarrow\Hom_{\mathcal{C}^\alpha}(Q^\alpha(\lambda),N^\alpha)\longrightarrow\Hom_{\mathcal{C}^\alpha}(Q^\alpha(\lambda),M^\alpha)\longrightarrow\Hom_{\mathcal{C}^\alpha}(Q^\alpha(\lambda),(M/N)^\alpha)\longrightarrow0.$$
for all $\lambda\in\Omega$ and $\alpha\in R^+$.
Then we have $$\mathbb{V}(M/N)(\lambda)=\mathbb{V}M(\lambda)/\mathbb{V}N(\lambda),$$ and
$$\mathbb{V}(M/N)(\lambda,\beta)=\mathbb{V}M(\lambda,\beta)/\mathbb{V}N(\lambda,\beta)$$ 
for any $\lambda\in\Omega$ and $\alpha\in R^+$.
Consider as $\hat{S}^\alpha$-module, we have $$(\mathbb{V}(M)(\lambda)\oplus\mathbb{V}(M)(\alpha\uparrow\lambda))/\mathbb{V}(N)(\lambda,\beta)\simeq S^\emptyset/S^\alpha\simeq(\mathbb{V}(N)(\lambda)\oplus\mathbb{V}(N)(\alpha\uparrow\lambda))/\mathbb{V}(N)(\lambda,\beta).$$
Then we have natural embeeding from
 $\mathbb{V}M(\lambda,\beta)/\mathbb{V}N(\lambda,\beta)$ into $\mathbb{V}(M)(\lambda)/\mathbb{V}(N)(\lambda)\oplus\mathbb{V}(M)(\alpha\uparrow\lambda)/\mathbb{V}(N)(\alpha\uparrow\lambda)$. It follows that $\mathbb{V}(M)/\mathbb{V}(N)\in\mathcal{K}$.
We conclude that $\mathbb{V}(M/N)\simeq\mathbb{V}M/\mathbb{V}N$ as desired.
\end{proof}

Let $\mathcal{Z}_\lambda=\mathbb{V}Z_{\hat{S}}(\lambda)$ for each $\lambda\in X$. This Lemma shows that 
If $M\in\mathcal{FC}_{\hat{S}}(\Omega)$ has a $Z$-filtration, then $\bbv(M)$ has a filtration with all factors of form $\mathcal{Z}_\lambda$ with $\lambda\in\Omega$. 
We call it a $\mathcal{Z}$-filtration.
Suppose $\mathcal{M}$ has a $\mathcal{Z}$-filtration. Then for any $\phi\in\Hom_{\mathcal{K}}(\mathcal{M},\mathcal{Z}_\mu)$, the rank of all nonzero $\phi_\lambda$ over $\hat{S}^\emptyset$ are equal to each other and we denote it as $\rank_\mathcal{K}\Hom_{\mathcal{K}}(\mathcal{M},\mathcal{Z}_\mu)$.
By Lemma \ref{lem proj over local ring}, we have projective object $Q_{\hat{S}}(\lambda)\in\mathcal{C}_{\hat{S}}$ such that $Q_{\hat{S}}(\lambda)\otimes k=Q_k(\lambda)$.

\begin{thm}
Keep notations as above, and let $\lambda,\mu\in\Omega$. Then the  multiplicity of $L_k(\mu)$ in baby Verma module $Z_k(\lambda)$ is equal to the multiplicity of $\mathcal{Z}_\mu$ in $\mathcal{Z}$-filtration of $\mathbb{V}(Q_{\hat{S}}(\lambda))$.
In particular, we have
$$[Z_k(\mu)\colon L_k(\lambda)]=\rank_{\hat{S}^\emptyset}\;\mathbb{V}Q_{\hat{S}}(\lambda)(\mu).$$
\end{thm}
\begin{proof}We have following equalities$\colon$
$$\begin{aligned}
[Z_k(\mu)\colon L_k(\lambda)]&=\dim_k\Hom_{\mathcal{C}_k}(Q_k(\lambda), Z_{k}(\mu))\\
&=\rank_{\hat{S}}\Hom_{\mathcal{C}_{\hat{S}}}(Q_{\hat{S}}(\lambda), Z_{\hat{S}}(\mu))\\
&=\rank_{\mathcal{K}}\Hom_{\mathcal{K}}(\mathbb{V}Q_{\hat{S}}(\lambda),\mathcal{Z}_\mu)\\
&=\rank_{\hat{S}^\emptyset}\Hom_{\mathcal{C}^\emptyset}(\mathbb{V}Q_{\hat{S}}(\lambda)(\mu),\mathcal{Z}_\mu(\mu))\\
&=\rank_{\hat{S}^\emptyset}(\mathbb{V}Q(\lambda)(\mu)).
\end{aligned}$$
The first equation is the property of the projective object in $\mathcal{C}_k$. The second equation comes from the fact $\Hom_{\mathcal{C}_{\hat{S}}}(Q_{\hat{S}}(\lambda), Z_{\hat{S}}(\mu))\simeq\Hom_{\mathcal{C}_{\hat{S}}}(Q_{\hat{S}}(\lambda), Z_{\hat{S}}(\mu))\otimes k$. The third equation comes by Proposition \ref{prop fully faithful}. The fourth equation is followed by the definition. The last equation is deduced from the fact $\mathcal{Z}_\mu(\mu)=\hat{S}^\emptyset$.
 \end{proof}

\end{document}